\documentclass[12pt]{amsart} 
\newfont{\sheaf}{eusm10 scaled\magstep1}

\newcommand{\C}{\ensuremath{\mathbb{C}}}
\newcommand{\R}{\ensuremath{\mathbb{R}}} 
\newcommand{\Z}{\ensuremath{\mathbb{Z}}} 
\newcommand{\Q}{\ensuremath{\mathbb{Q}}}
\newcommand{\M}{\ensuremath{\mathbb{M}}}
 
\newcommand{\hol}{\ensuremath{\mathcal{O}}}
\newcommand{\HH}{\ensuremath{\mathbb{H}}} 
\newcommand{\PP}{\ensuremath{\mathbb{P}}} 
\newcommand{\RR}{\ensuremath{\mathcal{R}}}
\newcommand{\BB}{\ensuremath{\mathcal{B}}}
\newcommand{\FF}{\ensuremath{\mathcal{F}}}
\newcommand{\HHH}{\ensuremath{\mathcal{H}}}
\newcommand{\V}{\ensuremath{\mathcal{V}}}

\newcommand{\Proof}{{\it Proof. }} 
\newcommand{\ra}{\ensuremath{\rightarrow}}

\newtheorem{teo}{Theorem}[section] 
\newtheorem{df}[teo]{Definition} 
\newtheorem{lem}[teo]{Lemma} 
\newtheorem{cor}[teo]{Corollary} 
\newtheorem{oss}[teo]{Remark} 
\newtheorem{prop}[teo]{Proposition}
\newtheorem{ex}[teo]{Example}

\def\eea{\end{eqnarray*}}
\def\bea{\begin{eqnarray*}}

\begin{document}

\title{Deformation in the large of some complex manifolds, I}

\author{Fabrizio Catanese\\
 Universit\"at Bayreuth
 }

\footnote{
The present research took place in the framework of the Schwerpunkt
"Globale Methode in der komplexen Geometrie", and of the EAGER EEC
Project.  
}
\date{October 13, 2002}
\maketitle

{\bf This article is dedicated to the memory of Fabio Bardelli}\\

\section{Introduction}

 The main  theme of the present note is the study of the deformations in the
large of  compact complex manifolds. Even when these manifolds are
K\"ahler, we shall study their deformations without imposing the K\"ahler
assumption.

Recall in fact that (cf. e.g. \cite{k-m71}) a small deformation of a
K\"ahler manifold is again K\"ahler. This is however (cf. \cite{hir62})
false for non small deformations.

In general, Kodaira defined two complex manifolds $X'$, $X$ to be {\bf directly
deformation equivalent} if there is a  proper holomorphic submersion
$ \pi : \Xi \ra
\Delta$ of a complex manifold $\Xi$ to the unit disk in the complex plane, such that
$X, X'$ occur as fibres of $\pi$. If we take the equivalence relation generated
by direct deformation equivalence, we obtain the relation of {\bf deformation
equivalence}, and we say that $X$ is a deformation of $X'$ in the large if
$X, X'$ are deformation equivalent.

These two notions extend the classical notions of irreducible, resp. connected
components of moduli spaces. However, outside of the realm of projective
manifolds, not so much is known about deformations in the large of compact complex
manifolds, since the usual deformation theory considers only the problem of
studying the small deformations. 

Just to give an idea of how limited our knowledge is, consider that only in
a quite recent
 paper (\cite {cat02}), of which this one is a continuation, we  gave a
positive answer to a basic question raised by Kodaira and Spencer (cf.
\cite {k-s58}, Problem 8, section 22, page 907 of volume II of Kodaira's
collected works), showing that any deformation in the large of a complex
torus is again a complex torus (in \cite{k-s58} only the case
n=2 was solved, cf. Theorem 20.2).

The usual strategy to determine the deformation equivalence class of a
compact complex manifold $X$ is to construct as big a family $\Xi \ra \BB$
as possible, where $\BB$ is a connected analytic space, and
then try to prove
\begin{itemize}
\item
The family is versal, i.e., for each fibre $X_0$, we get a local surjection
$(\BB, 0) \ra Def (X_0)$ onto the Kuranishi family of $X_0$.
\item
Given any 1-parameter family $ \pi : \Xi' \ra
\Delta$ with the property that there exists a sequence $t_{\nu} \ra 0$
such that $X_{t_{\nu}}$,
$\forall \nu$, occurs as a fibre in the family
$\Xi$, then the same  property is also enjoyed by $X_0$ .
\end{itemize}

In this paper we shall  begin to go further, considering manifolds which are
torus fibrations, and since we will make  extensive use of  results and
techniques from (\cite {cat02}), in the second section we will reproduce and
extend  some of those for the reader's and ours benefit.
In particular one can find a complete proof of

{\bf Theorem 2.1 }: 
{\it Every deformation of a complex torus of dimension $n$ is a 
complex torus of dimension $n$.}\footnote{Marco Brunella pointed out that
a weaker result in this direction had been obtained by A. Andreotti and
W.Stoll in the paper "Extension of holomorphic maps", Annals of Math. 72,
2 (1960), 312-349.}

In section 3 we shall consider the general situation of a 
holomorphic fibration $ f : X \ra Y$ with all the fibres complex tori 
(equivalently, by theorem 2.1, with all the fibres smooth and with one fibre
isomorphic
to a complex torus).

We shall set 
up a standard  notation, and we shall give criteria for $f$ to be a holomorphic torus
bundle, respectively a principal holomorphic torus bundle. The main result is
proposition 3.4, describing explicitly all principal holomorphic torus
bundles over curves.

{\bf Proposition  3.4 }:
{\it Any principal holomorphic torus bundle $X$ over a curve $Y$ is a quotient
$X = L / N$ of a suitable holomorphic $(\C^*)^d$-bundle $L$ over $Y$ by a suitable
discrete cocompact subgroup $N$ of $(\C^*)^d$.}

Section 4 treats first the general problem of  determining the small deformations
of torus bundles, and then the main {\bf Theorem 4.4 } asserts that when the base is
a curve of genus $\geq 2$, then all the small deformations of a  principal
holomorphic torus bundle are again holomorphic torus bundles.

The short Section 5  is devoted to showing that the families 
 of   principal holomorphic torus bundles on curves constructed
in section 4 yield all the deformations in the large of such manifolds.

{\bf Theorem 5.1 }: 
{\it 
A  deformation in the large of a holomorphic principal torus bundle over
a curve $C$ of genus $g\geq 2$ with fibre a complex torus $T$ of
dimension $d$ is again a holomorphic principal torus  bundle
over a curve $C'$ of genus $g$.}

Section 6 is instead quite long and develops a classification theory for principal
holomorphic torus bundles  over tori. This
theory bears close similarities with the theory of line bundles over tori, as
follows. Namely, assume that we have such a torus bundle  over $Y = V/\Gamma$
and with fibre $T = U / \Lambda$.

Then the homomorphism of fundamental groups $ \pi_1(X) \ra \pi_1(Y) = \Gamma$
is a central extension with kernel $\Lambda$ and is completely classified
by a bilinear alternating form $ A : \Gamma \times \Gamma \ra \Lambda.$

We have a vector analogue of the Riemann bilinear relation, since
viewing $A$ as a real element of 
 $$\Lambda^2(\Gamma 
\otimes \R )^{\vee} \otimes (\Lambda \otimes \R )=  \Lambda^2(V \oplus
\bar{V})^{\vee} \otimes (U \oplus \bar{U}),$$
 its component in $  \Lambda^2(
\bar{V})^{\vee} \otimes (U )$ is zero.

The Riemann bilinear relation enables us to construct a
universal family of such bundles, for each  given choice of the extension class 
$A$. Namely, we take first the possible subspaces $ U \subset (\Gamma 
\otimes \R )$ , $V \subset (\Lambda \otimes \R )$ such that the Riemann bilinear
relation holds, and in this way we obtain a family which we call the
standard Appell Humbert family.

From this family we obtain the so called complete Appell Humbert family,
which contains all such principal holomorphic bundles with given form
$A$, as stated by 

{\bf  Theorem 6.8}
{ \it Any  holomorphic principal torus bundle with extension class isomorphic
to $\epsilon  \in H^2 (\Gamma, \Lambda)$ occurs in the complete
Appell-Humbert family  $\mathcal T' \mathcal B_A$.}

The similarity
with the theory of line bundles on tori, namely with the Theorem of
Appell-Humbert, is that we are able to  explicitly write the classifying 
group cocycle  as a linear function easily determined by the extension
class and  by the complex structures
$V$  on the base, respectively $U$ on the fibre.

We show later the advantages of having such a realization of these bundles
via explicit cocycles, namely how one can explicitly calculate several
holomorphic invariants of such bundles using the bilinear algebra data of
the extension class, we have for instance

{\bf  Theorem 6.10}
{ \it 
The cokernel of $ 0 \ra H^0 (\Omega^1_Y )\ra H^0 (\Omega^1_X )$ is the
subspace of ${U}^{\vee}$ which annihilates the image of the Hermitian
part of $A$, i.e.,  of the component $B"$ in 
$[( V \otimes \bar{V})^{\vee}  \otimes (U )].$

It follows in particular that $X$ is parallelizable if and only if the
Hermitian part of $A$ is zero.}

{\bf  Corollary 6.11}
{\it The space $H^0 ( d\hol_X)$ of closed holomorphic 1-forms on $X$ contains 
the pull- back of $H^0 ( \Omega^1_Y)$ with cokernel the subspace $U^*$
of $U^{\vee}$ which annihilates the image of $A$, i.e.,  $U^* = \{ \beta | 
\beta
\circ A (z,\gamma)= 0 \
\forall
\gamma, \ \forall z\} .$}

{\bf  Theorem 6.12}
{ \it 
The cokernel of $ 0 \ra H^1 (\hol_Y )\ra H^1 (\hol_X )$ is the subspace
of $\bar{U}^{\vee}$ which annihilates the image of the anti-complex component $B'$ 
 of $A$, i.e., of the conjugate of the component in 
$[( \Lambda^2 V )^{\vee}  \otimes (U )].$ }

As the reader might have noticed, we give necessary and sufficient conditions for the
parallelizability of
$X$, showing thus (a fact pointed out  by Nakamura already many years ago, cf.
\cite{nak75}) easily that this notion is not deformation invariant. 

Via similar calculations, we  determine  (Proposition 6.16) the
cohomology algebra  $ \oplus _i H^i ( \hol_X)$ of the structure sheaf,
and ( Corollary 6.15) the cohomology groups of the Tangent sheaf of $X$.

With these calculations we are able to show  (Theorem 6.17) that the complete
Appell Humbert family is versal for the manifolds $X$ corresponding to smooth
points of the base of the family, if suitable assumptions hold for $A$, and in
the case where the fibre dimension is $1$.

For lack of time we defer to the future the investigation of the small deformations
in the general case, which is a necessary step for the investigation of the
deformations in the large.


We finish the section by considering in detail the classical example
of the Iwasawa 3-fold, and of its small deformations.

Finally, in section 7, we recall the definition (cf. \cite{cat02}) of
Blanchard-Calabi torus fibrations, and  using the result of
section 5 we obtain as a consequence that the space of complex structures
on the differentiable manifold underlying the product of a curve with a
complex torus of dimension 2 has several distinct deformation types, namely we have
the following

{\bf Corollary 7.8}

{\it The space of complex structures on the product of a curve $C$
of genus $g \geq 2$ with a four
dimensional real torus contains manifolds which are not deformation
equivalent to each other, namely, Blanchard-Calabi 3-folds which are not
K\"ahler and holomorphic principal bundles in the family $\FF _{g,0}$.}

\section{Tori}

 It is well known that complex tori are parametrized by a connected
family (with smooth base space), inducing  all the small deformations (cf.
\cite{k-m71}).

In fact, we have a family parametrized by an open set $\mathcal T_n$ of the complex
Grassmann Manifold $Gr(n,2n)$, image of the open set of matrices

$\{ \Omega \in Mat(2n,n; \C)| det  (\Omega \overline {\Omega}) > 0 \}.$

We recall this well known fact: if we consider $\Z^{2n}$ as a fixed lattice, to each
matrix $ \Omega $ as above we associate the subspace $ U = ( \Omega ) (\C^{n})$, so
that 
$ U \in Gr(n,2n)$ and $\Z^{2n} \otimes \C \cong U \oplus \bar{U}.$

Finally, to $ \Omega $ we associate the torus $U / p_U (\Z^{2n})$, $p_U : U \oplus
\bar{U} \ra U$ being the projection onto the first addendum.

The above family will be the main object we shall have in mind.

We also want to recall from (\cite {cat02})  the following first result:

\begin{teo}
Every deformation of a complex torus of dimension $n$ is a 
complex torus of dimension $n$.
\end{teo}

In fact this proof  uses  several useful lemmas (some of them well known)
 which will  be here  slightly generalized, for further use  in the
sequel.  The first two are due to Kodaira (for the first, cf. \cite{ko64},
Theorem 2, page 1392 of collected works, vol. III)

\begin{lem}
On a compact complex manifold $X$ one has an injection

$$ H^0( d \hol_X ) \oplus  \overline{H^0( d {\hol_X })}   
\rightarrow  H^1_{DR} (X,  \C).$$
\end{lem}

PROOF. Suffices to show that the map $ H^0( d \hol_X )   
\rightarrow  H^1_{DR} (X,  \R)$, sending 
$$ \omega  \rightarrow \omega  + \bar{\omega } $$is an injection.
Else, there is a function $f$ with $df = \omega  + \bar{\omega }$,
whence $ \partial f = \omega $ and therefore $ \bar{\partial} 
\partial f =
d (\omega) = 0$. Thus $f$ is pluriharmonic, hence constant by the maximum
principle. Follows that $\omega = 0$.
\hfill Q.E.D. \\

Recall (cf.  \cite {bla56} and \cite {ue75}) that for a compact complex manifold
$X$, the Albanese  Variety
$Alb(X)$ is the quotient of the complex dual vector space of $H^0( d \hol_X 
)$ by the minimal  closed complex Lie subgroup containing the image of
$H_1(X, \Z)$.

The Albanese map
$\alpha_X : X \rightarrow Alb(X)$ is given as usual by fixing a
base point $x_0$, and defining $\alpha_X (x) $ as the  class in the
quotient of the functional given by integration on any  path connecting $x_0$ with
$x$.

One says the the Albanese Variety is {\bf good} if the image $H_1(X)$ of
$H_1(X, \Z)$ is discrete in $H^0( d \hol_X  )$, and {\bf very good } if it is a
lattice (a discrete subgroup of maximal rank). Moreover, the {\bf Albanese
dimension of $X$ } is defined as the dimension of the image of the Albanese map.

With this terminology, we can state an important consequence of an inequality 
due to Kodaira

\begin{lem}
On a compact complex manifold $X$ one has an injection

$$  \overline{H^0( d {\hol_X })}   
\rightarrow  H^1 (X,  \hol_X).$$
In particular, if $b_1(X)$ is the first Betti number of $X$, we have the
inequalities
$$2 h^1 (X,  \hol_X) \geq h^1 (X,  \hol_X) + h^0( d {\hol_X })
\geq b_1(X).$$
If both equalities hold, then $X$ has a very good Albanese variety, of
dimension
$h^1 (X,  \hol_X) = h^0( d {\hol_X }) = \frac{1}{2} 
 \ b_1(X)$.
\end{lem}

PROOF. We claim that the map $ \overline {H^0( d \hol_X )}   
\rightarrow  H^1 (X,  \hol_X) = H^1_{\bar \partial} (X,  \hol_X)$,  
is  injective.
Else, we have $  \bar{\omega } $  with $ \partial \bar{\omega }=0 $
and $ \bar \partial \bar{\omega }=0 $ such that
there is a function $f$ with $\bar \partial f =   \bar{\omega }$,
whence  $ \bar{\partial} 
\partial f =
- d (\omega) = 0$. Thus we conclude as in the preceding lemma: $f$ is
pluriharmonic, whence constant, thus $\omega = 0$.

The second assertion follows easily from the exact cohomology sequence
$$  H^0( d \hol_{X} ) \ra H^1 (X,  \C) \ra 
H^1(X,  \hol_{X} )  $$
and from the first.

Finally, if equality holds it follows that the injection 
$$ H^0( d \hol_X ) \oplus  \overline{H^0( d {\hol_X })}   
\rightarrow  H^1_{DR} (X,  \C)$$
is an isomorphism, whence $H_1 (X,  \Z)/ (Torsion )$ maps isomorphically to
a lattice  $ H_1(X)$  in the dual
vector space of $ H^0( d \hol_X )$, thus $Alb(X) =  H^0( d \hol_X
)^{\vee}/(H_1 (X)).$
\hfill Q.E.D. \\

\begin{lem}
Assume that $\{X_t\}_{t \in \Delta}$ is a 1-parameter family of compact
complex manifolds over the 1-dimensional disk, such that there is a
sequence $ t_{\nu} \rightarrow 0$ with $X_{t_{\nu}}$ 
satisfying the weak 1-Hodge  decomposition 
$$ H^1_{DR} (X_0,  \C) = H^0( d \hol_{X_0} ) \oplus  
\overline {H^0({d \hol_{X_0}} )}.  $$

Then the weak 1-Hodge  decomposition
holds also on the central fibre $X_0$.
\end{lem}

PROOF. We have $ f: \Xi \rightarrow \Delta  $ which is proper and smooth,
and  $ f_* (\Omega^1_{\Xi |\Delta}) $  is torsion free, whence ($\Delta$ is
smooth of dimension $1$) it is locally free of rank $ h \geq q: =
(1/2) \ b_1(X_0).$

In fact, the Base change Theorem asserts that there is ( cf. \cite{mum70}, II 5, and
\cite{gr84}) a complex of Vector Bundles on
$\Delta$,

$$  (*)\  \  E^0 \rightarrow E^1 \rightarrow E^2 \rightarrow ...E^n\ 
s.t. $$ 
$1) R^if_* (\Omega^1_{\Xi |\Delta}) $ is the i-th cohomology group of (*),
whereas 

$ 2)H^i (X_t, \Omega^1_{X_t}) $ is the i-th cohomology group of
$(*) \otimes \C_t .$

We may shrink the disk $\Delta$ so that the rank of $E_i \ra E_{i+1}$ is constant
for $ t \neq 0$, and thus as a consequence, there is an isomorphism
between $H^0 (X_t, \Omega^1_{X_t}) $ and the stalk $ f_* (\Omega^1_{\Xi |\Delta})
 \otimes \C_t$. 
 
 We first see what happens instead on the central fibre, where the
 space $H^0 (X_0, \Omega^1_{X_0}) $ can have a higher dimension

CLAIM 1. 

There are holomorphic 1-forms $\omega_1(t), ...\omega_h(t)$
defined  in the inverse image $f^{-1} (U_0)$ of a neighbourhood 
$(U_0)$ of $0$ , and such that their restriction to $X_t$ are 
linearly independent  $\forall  t \in U_0$.

PROOF OF  CLAIM 1: assume that  $\omega_1(t), ...\omega_h(t)$ generate
the direct image sheaf  $ f_* (\Omega^1_{\Xi |\Delta}) $, but
$\omega_1(0), ...\omega_h(0)$ be linearly dependent. Then, w.l.o.g. we
may assume
$\omega_1(0)\equiv 0$, i.e. there is a maximal $m$ such that
$\hat{\omega}_1(t) : = \omega_1(t)/t^m$ is holomorphic. Then, since 
$\hat{\omega}_1(t) $ is a section of  $ f_* (\Omega^1_{\Xi |\Delta})
$, there are holomorphic functions $\alpha_i$ such that 
$\hat{\omega}_1(t)  = \Sigma_{i=1,..h} \alpha_i(t)\omega_1(t) $,
whence it follows that $\hat{\omega}_1(t) ( 1- t^m \alpha_1(t)) =
\Sigma_{i=2,..h} \alpha_i(t)\omega_1(t).$

This however contradicts the fact that $ f_* (\Omega^1_{\Xi |\Delta}) $ 
is locally free of rank $ h$.

CLAIM 2:

$H^0( d \hol_{X_0}) $ has dimension $\geq q$.

PROOF OF CLAIM 2.

Let $d_v$ be the vertical part of exterior differentiation, i.e., the
composition of $d :\Omega^1_{\Xi} \ra \Omega^2_{\Xi} $ with the projection $ 
(\Omega^2_{\Xi} \ra \Omega^2_{\Xi |\Delta}) $. 

It clearly factors through an $f^* \hol_{\Delta}$-linear map
$d'_v : \Omega^1_{\Xi |\Delta} \ra \Omega^2_{\Xi |\Delta}.$

Taking direct images, we get a homomorphism of coherent sheaves
$f_* (d'_v ): f_* (\Omega^1_{\Xi |\Delta}) \ra f_* (\Omega^2_{\Xi |\Delta})$
whose kernel will be denoted by $\mathcal H$.

By our assumption, for each $t_{\nu} $ we have  a $q$-dimensional subspace
$H_{t_{\nu}} := H^0( d \hol_{X_{t_{\nu}}})$ of $H^0 (X_{t_{\nu}},
\Omega^1_{X_{t_{\nu}}}) $. 

Therefore, as $ t_{\nu} \ra
0$, we have a limit in the Grassmann manifold  $Grass (q,h)$ of the family
$ \{ H_{t_{\nu}} \}$,
whence  a subspace $H_0$ of
dimension $q$. Taking suitable bases, $\omega_i(0)$ of $H_0$,
resp. $\omega_i(t_{\nu})$  of $H_{t_{\nu}}$ we see that, since
$\omega_i(t_{\nu})$ is $d_v$-closed,  by continuity it follows that  also
$\omega_i(0)
\in H^0 ( d\hol_{X_{0}}).$

END OF THE PROOF:

It follows from  lemma 2.2 that, being $b_1 = 2
q$, $ H^1_{DR} (X_0,  \C) = H^0( d \hol_{X_0}  ) \oplus  
H^0( d \bar{\hol}_{X_0} )$.

Actually, by the base change theorem follows that $\mathcal H$ is free of rank
$q$ and it enjoys the base change property that its stalk 
$\mathcal H \otimes \C_t $ corresponds to the subspace $ H^0( d
\hol_{X_t}).$

\hfill Q.E.D. \\

\begin{cor}
Assume that $\{X_t\}_{t \in \Delta}$ is a 1-parameter family of compact
complex manifolds over the 1-dimensional disk, such that there is a
sequence $ t_{\nu} \rightarrow 0$ with $X_{t_{\nu}}$ satisfying the weak
1- Hodge property, and moreover such that $X_{t_{\nu}}$  has  Albanese
dimension $=a$.

Then the central fibre $X_0$ has a very good Albanese Variety, and has
also  Albanese dimension $=a$.
\end{cor} 

PROOF. We use the fact (cf. \cite{cat91}) that, when the Albanese Variety is
good, then  the Albanese dimension of
$X$ is equal to $ max \{  i|\ \Lambda^i H^0( d \hol_X  ) \otimes 
\Lambda^i H^0( d \bar{\hol}_X  ) \rightarrow H^{2i}_{DR}
(X, \C)$ has non zero image $\}.$

If the weak 1-Hodge decomposition holds for $X$, then the Albanese
dimension of $X$ equals $(1/2) \   max \{  j| \Lambda^j H^1 (X, 
\C) $ has non zero image in   $ H^j (X, 
\C)\}.$  But this number is clearly invariant by
homeomorphisms.

Finally, the Albanese Variety for $X_0$ is very good since  the weak
1-Hodge decomposition holds for $X_{t_{\nu}}$ whence also for
$X_0$ by the previous lemma.

\hfill Q.E.D. \\ 

\begin{oss}
As observed in (\cite{cat95}, 1.9), if a complex manifold $X$ has a
generically finite map to a K\"ahler manifold, then $X$ is bimeromorphic to
a  K\"ahler manifold. This applies in particular to the Albanese map.
\end{oss}

Theorem 1.1 follows from the following statement

\begin{teo}
Let $X_0$ be a
compact complex manifold such that its Kuranishi family of deformations
$\pi : \Xi
\rightarrow
\mathcal{B}$ enjoys the property that the set $\mathcal{B}(torus) : = \{b|
X_b $ is isomorphic to a complex torus$\}$ has $0$ as a limit point.

Then $X_0$ is a complex torus.
\end{teo}

We will use the following "folklore"

\begin{lem}
Let $Y$ be a connected complex analytic space, and $Z$ an open set of $Y$
such that $Z$ is closed for holomorphic 1-parameter limits ( i.e., given
any holomorphic map of the 1-disk $f : \Delta \rightarrow Y$, if there is
a sequence $t_{\nu} \rightarrow 0$ with $f(t_{\nu} ) \in Z$, then also
$f(0) \in Z$). Then $Z=Y$.
\end{lem}

PROOF OF THE LEMMA. 

By choosing an appropriate stratification of $Y$ by smooth manifolds, it
suffices to show that the statement holds for $Y$ a connected manifold.

Since it suffices to show that $Z$ is closed,  let $P$ be a
point in the closure of $Z$, and let us take coordinates such that a
neighbourhood of $P$ corresponds to a compact polycylinder $H$ in $\C^n$.

Given a point in $Z$, let $H'$ be a maximal coordinate polycylinder
contained in $Z$. We claim that $H'$ must contain $H$, else, by the
holomorphic 1-parameter limit property, the boundary of $H'$ is contained in
$Z$, and since $Z$ is open, by compactness we find a bigger polycylinder
contained in $Z$, a contradiction which proves the claim.

\hfill Q.E.D. for the Lemma.\\ 

PROOF OF  THEOREM 2.7.

It suffices to consider a 1-parameter family ($\mathcal{B} = \Delta$)
whence we may assume w.l.o.g. (cf. the proof of Lemma 2.4) that the weak
1-Hodge decomposition holds for each
$t \in \Delta$.

By integration of the holomorphic 1-forms on the fibres ( which are
closed for $t_{\nu}$ and for $0$), we get a family of Albanese maps
$\alpha_t : X_t
\rightarrow J_t$, which fit together in a relative map $a : \Xi
\rightarrow J$ over $\Delta$ ( $J_t $ is the complex torus
$(H^0(d \hol_{X_t}))^{\vee}/ H_1(X_t, \Z)$).

Apply once more  vertical exterior differentiation to the forms 
$\omega_i(t)$: $ d_v (\omega_i(t)$  vanishes identically on $X_0$ and on
$X_{t_{\nu}}$, whence it vanishes identically in a neighbourhood of
$X_0$, and therefore these forms $\omega_i(t)$ are closed for each $t$.

Therefore our map $a : \Xi \rightarrow J$ is defined everywhere and it is an
isomorphism for $t = t_{\nu}.$

Whence, for each $t$, $\alpha_t$ is surjective and has degree 1.

To show that $\alpha_t$ is an isomorphism for each $t$ it suffices
therefore to show that $a$ is finite. 

Assume the contrary: then there is a ramification divisor $R$ of $a$,
which is exceptional (i.e., if $B = a(R)$, then $dim B < dim R$).

By our hypothesis $\alpha_t$ is an
isomorphism for $t = t_{\nu}$, thus $R$ is contained in a union of fibres,
and since it has the same dimension, it is a finite union of fibres.
But if $R$ is not empty, we reach a contradiction, since then there are
some $t$'s such  that $\alpha_t$ is not surjective.

\hfill Q.E.D. \\ 

We have also a more abstract result 

\begin{prop}
Assume that $X$ has the same integral cohomology algebra of a complex
torus and that   $ H^0( d \hol_X  )$ has dimension equal to $n =dim(X)$.
Then $X$ is a complex torus.
\end{prop}

PROOF.

Since $b_1(X) = 2n$, it follows from 1.2 that the weak 1-Hodge
decomposition holds for $X$ and that the Albanese Variety of $X$ is very
good.

That is, we have the Albanese map $\alpha_X : X \rightarrow J$, where 
$J$ is the complex torus $J= Alb(X)$. We want to show that the Albanese
map is an isomorphism. It is a morphism of degree $1$, since $\alpha_X $
induces an isomorphism between the respective fundamental classes of
$H^{2n}(X,\Z) \cong H^{2n}(J,\Z)$.

There remains to show that the bimeromorphic morphism $\alpha_X $ is finite.

To this purpose, let $R$ be the ramification divisor of $\alpha_X : X
\rightarrow J$, and $B$ its branch locus, which has codimension at least
$2$. By means of blowing ups of $J$ with non singular centres we can
dominate
$X$ by a K\"ahler manifold $g: Z \rightarrow X $ (cf. \cite{cat95} 1.8, 1.9).

Let $W$ be a fibre of $\alpha_X$ of positive dimension such that
$g^{-1} W$ is isomorphic to $W$. Since $Z$ is K\"ahler, $g^{-1} W$ is not
homologically trivial, whence we find a differentiable submanifold $Y$ of
complementary dimension which has a positive intersection number with it.
But then, by the projection formula, the image $g_* Y$ has positive
intersection with
$W$, whence $W$ is also not homologically trivial. However, the image of
the class of $W$ is $0$ on $J$, contradicting that $\alpha_X$ induces an
isomorphism of cohomology (whence also homology) groups.

\hfill Q.E.D. \\ 

\begin{oss}
The second condition holds true as soon as the complex dimension $n$ is
at most $2$. For $n=1$ this is well known, for $n=2$ this is also known,
and due to Kodaira (\cite{ko64}): since for $n=2$ the holomorphic 1-forms
are closed, and moreover  $ h^0( d \hol_X  )$ is at least $[ (1/2)
b_1(X)]$.

For $n\geq 3$, the real dimension of $X$ is greater than $5$, whence, by
the s-cobordism theorem (\cite{maz63}), the assumption that $X$ be
homeomorphic to a complex torus is equivalent to the assumption that  $X$ be
diffeomorphic to a complex torus. 

Andr\'e Blanchard (\cite{bla53}) constructed in the early 50's an example of
a non K\"ahler complex structure on the product of a rational curve with a
two dimensional complex torus. In particular his construction (cf.
\cite{somm75}), was rediscovered by Sommese, with a more clear and
 more general presentation, who pointed out that in this way one would
produce exotic complex structures on complex tori. 

There remains however open the following

 QUESTION: let $X$ be a  compact complex manifold of dimension $n \geq
3$  and with trivial canonical bundle such that $X$ is diffeomorphic to a
complex torus: is then $X$ a complex torus ?
\end{oss}

The main problem here is to show the existence of holomorphic 1-forms, so
it may well happen that also this question has a negative answer.

\section{Generalities on holomorphic torus bundles}

Throughout the rest of the paper, our set up will be the
following: we have a holomorphic submersion between compact
complex manifolds

$$ f : X \ra Y  ,$$
such that one fibre $F$  (whence all the fibres, by theorem 2.1)
is a complex torus.

We shall denote this situation by saying that $f$ is 
differentiably a torus bundle.

We let $n= dim X$, $m = dim Y$, $d=dim F= n -m$.

The case where $d=1$ is very special because the moduli space
for
$1$-dimensional complex tori exists and is isomorphic to $\C$.
Whence it follows that in this situation (unlike the case $d
\geq 2$) $f$ is a holomorphic fibre bundle.

In any dimension, we have (e.g. by a much more general theorem
of Grauert and Fischer ( cf. \cite {fg65})) that $f$ is a
holomorphic bundle if and only if all the smooth fibres are 
biholomorphic.

For any differentiable torus bundle we have a local system on
$Y$,
$$ (1) \ \HH : = \RR^1 f_* \Z_X $$
and, by the exponential sequence $0 \ra \Z_X \ra \hol_X \ra
\hol^*_X \ra 0$, if we define 
$$ (2) \ \V^{\vee}: = \RR^1 f_* \hol_X  ,$$
then we get another exact sequence
$$ (3) \ 0 \ra W^{\vee} \ra \HH \otimes  \hol_Y \ra \V^{\vee}
\ra 0 $$
of holomorphic vector bundles, where
$$ (4) \ W^{\vee}: =  f_* \Omega^1_{X|Y}  .$$

Here, $\HH$ yields an $\R$-basis of $ \V^{\vee}$ at each point, and the
same holds for $ \HH^{\vee} : = { \it Hom}_{\Z} (\HH, \Z)$ via the exact
sequence
$$ (5) \ 0 \ra \V \ra \HH^{\vee} \otimes  \hol_Y \ra W
\ra 0 $$
so that one obtains another  differentiable bundle of complex tori,
the so called Jacobian of $X$, which is a bundle of Lie groups
$$ Jac (X) := W / \HH^{\vee}  . $$ 
By Hodge Symmetry the  real flat bundle $\HH^{\vee} \otimes  \hol_Y $ splits
as a direct sum $ \V \oplus \bar{\V}$, and $ W \equiv \bar{\V}$ as complex vector
bundles.
We finally have the cotangent sheaves exact sequence
$$(6) \  0 \ra f^* \Omega^1_Y \ra  \Omega^1_X \ra f^* W^{\vee} \ra 0$$ 
by which it follows
$$ (7) \ K_X \equiv f^* (  K_Y + det W^{\vee} ).$$
In general, the derived direct image cohomology sequence of (6) 
$$(8) \  0 \ra  \Omega^1_Y \ra f_* \Omega^1_X \ra f_* \Omega^1_{X|Y}
\ra  \Omega^1_Y \otimes \RR^1 f_* \hol_X \ra ..$$
is such that the coboundary map is given by the Kodaira-Spencer class
in $H^0 (\Omega^1_Y \otimes \RR^1 f_* \HHH om(\Omega^1_{X|Y},\hol_X
))$ which vanishes exactly when $f$ is a holomorphic fibre bundle.
In this case we have then an exact sequence
$$(8) \  0 \ra  \Omega^1_Y \ra f_* \Omega^1_X \ra f_* \Omega^1_{X|Y}
\ra 0 .$$

\begin{oss}
An immediate corollary of (6) is that $X$ is complex parallelizable, i.e.,
$\Omega^1_X$ is trivial,   only if $Y$ is parallelisable and $W$ is trivial
on $Y$. If $\Omega^1_Y$ and $W$ are trivial, then $X$ is parallelizable if
furthermore the extension class of $(6)$,  is trivial. This extension
class lies  in $H^1 ( f^* (\Omega^1_Y \otimes W )).$ If moreover we have
a holomorphic bundle, the class lies in $H^1 ( \Omega^1_Y \otimes W)$
and $X$ is complex parallelizable if and  only if $Y$ is parallelisable,
$W$ is trivial on $Y$, and (8) splits.
\end{oss}

Assume now that we have a holomorphic torus fibre bundle, thus we have the
exact sequence 
$$(8') \  0 \ra  \Omega^1_Y \ra f_* \Omega^1_X \ra W^{\vee} \ra 0 .$$

We have  the well known

\begin{prop}
A differentiable bundle of $d$-dimensional tori is a principal
holomorphic bundle if and only if we have an exact sequence
$$(8') \  0 \ra  \Omega^1_Y \ra f_* \Omega^1_X \ra \hol_Y^d \ra 0 .$$
\end{prop}

\Proof
$f$ is a holomorphic bundle if and only if the Kodaira Spencer class
is identically zero, i.e., (8) is exact.
In general (cf. \cite{bpv84}) if $T$ is a complex torus, we have an
exact sequence of complex Lie groups
$$ 0 \ra T \ra Aut(T) \ra \M \ra 1 $$
where $\M$ is discrete.
Taking sheaves of germs of holomorphic maps with source $Y$ we get
$$ 0 \ra \HHH(T)_Y \ra \HHH (Aut(T))_Y \ra \M \ra 1 $$
and we know that holomorphic bundles with base $Y$ and fibre $T$ are
classified by the cohomology group $H^1(Y, \HHH (Aut(T))_Y)$.
The exact sequence 
$$ 0 \ra H^1(Y, \HHH (T)_Y) \ra H^1(Y, \HHH (Aut(T)))_Y \ra H^1(Y,
\M) $$
determines when a holomorphic bundle is a principal holomorphic bundle.
In this case the cocycles  are in $H^1(Y, \HHH (T)_Y)$, i.e., they have values in the
translation group, whence $W$ is a trivial bundle.

Conversely, if $W$ is  trivial, we may first choose local
coordinates
$(y)$ on a small neighbourhood $U \subset Y$ and local coordinates $ (x)
= (u',y)$ on $X$ with
$ f(x) = y$, then we may choose a basis $ w_1(y), .. w_d(y)$ of $
W^{\vee}$, lift these to local holomorphic
$1$-forms on $f^{-1} (U)$, $ \omega_1(y), .. \omega_d(y)$ and then take 
 the linear coordinates
$u_i :=  \int^{x}_{(0,y) } \omega_i (u',y)$.

On the universal cover of $f^{-1}(U)$,  we get functions $(y_1,..
y_m),(u_1, .. u_d)$ whose differentials give a basis of $\Omega^1_X$. 
Moreover, if we go to another open set $V \subset Y$, the new linear
coordinates $(v_1, .. v_d)$ are such that $v_i - u_i$ is a function
of $y \in U \cap V$.
\qed

In the case of a principal holomorphic bundle it is useful to
write $T = \C^d / \Lambda$, where $\Lambda \cong \Z^{2d}$, thus we have
the exact sequence
$$(9) \   \ra H^0(\HHH(T)_Y) \ra H^1(Y, \Lambda ) \ra  H^1 (Y,
\hol_Y^d) \ra H^1(\HHH(T)_Y) \ra ^c {\ra} H^2(Y, \Lambda ) .$$
The significance of the homomorphism $c$ is readily offered by
topological considerations.

In fact, the homotopy exact sequence of a bundle $f$
$$ (10) \pi_2(Y) \ra \pi_1(F) \ra \pi_1(X) \ra \pi_1(Y) \ra 1$$
determines an extension 
$$ (10') 1 \ra \pi_1(F) \ra \Pi \ra \pi_1(Y) \ra 1$$
where the action of $\pi_1(Y)$ on $\pi_1(F)$ by conjugation is 
precisely the monodromy automorphism.
For a principal torus bundle the monodromy is trivial, indeed
$\HH \otimes \hol_Y$ is a trivial differentiable bundle,
thus the extension (10') is central, and is therefore classified
by ($ \Lambda = \pi_1(T)$, and $F = T$) a cohomology class
$\epsilon \in H^2(Y, \Lambda )$.

The following situation is very interesting: let us consider
 a Stein manifold $Y'$  whose cohomology group $H^2(Y', \Z)= 0$,
 and a holomorphic map $ h : Y' \ra Y$. Then we may replace
our principal bundle with compact base $Y$,  with its pull-back via $h$.
It follows from the previous considerations that the pull back
is a product $ Y' \times T$ (since $ \Lambda \cong \Z^d$,
and in view of (9) ). Actually, this holds more generally if
the pull back to $Y'$ of $\epsilon$ is trivial.

Noteworthy special cases may be :
\begin{itemize}
\item
$Y'$ is the universal cover of $Y$, e.g. if $Y$ is a curve or
a complex torus.
\item
$Y'$ is an open set of $Y$, for instance, if $Y$ is a curve, any 
open set $ \neq Y$ satisfies our hypotheses.
\item
$Y'$ is $(\C^*)^m$ and the pull back of the class $\epsilon$ is trivial.
\end{itemize}

Nevertheless, the case where $Y$ is a curve can be more easily
described via the following construction (cf. e.g. \cite{bpv84} pp.
143-149)

\begin{ex}
Let $Y$ be a complex manifold and let $L$ be a principal
holomorphic $(\C^*)^d$-bundle over $Y$, thus classified by a cohomology
class $\xi \in H^1(Y, \hol^*_Y)^d$.
Let us consider any complex torus $T = \C^d / \Lambda$, and let us
realize it as a quotient $(\C^*)^d / N$.
Then we can form the quotient $L / N$, which is a holomorphic fibre
bundle with fibre $\cong T$.
\end{ex}

\begin{prop}
Any principal holomorphic torus bundle $X$ over a curve $Y$ is a quotient
$X = L / N$ of a suitable holomorphic $(\C^*)^d$-bundle $L$ over $Y$.
\end{prop}

\Proof
It suffices (cf. loc. cit.) to show that there is a primitive embedding
$ i: \Z^d \ra \Lambda$ and a class $\eta$ in $H^2 (Y, \Z^d)$ such that
the class
$\epsilon$ corresponding to $X$ equals $i_* (\eta)$.
Consider however that $H^2 (Y, \Lambda) \cong H^2 (Y, \Z) \otimes
\Lambda$: thus there is a primitive embedding of $\Z $ in $\Lambda$
with the property that the fundamental class of $Y$ maps to $\epsilon$.
It suffices to extend the primitive embedding of $\Z$ to one of $\Z^d$
into $\Lambda$, and define $\eta$ as the image of the fundamental class
for the homomorphism induced by the embedding $ \Z \subset \Z^d$.
\qed

\section{Small deformations of  torus bundles}

Let us  consider again briefly the case of a general  differentiable torus
bundle.

Consider the exact sequence 
$$ (12)\ 0 \rightarrow f^* (W) \rightarrow  \Theta_X   \rightarrow  f^*
(\Theta_Y)  \rightarrow 0 $$
and the derived direct image sequence  
$$(13) 0 \rightarrow (W) \rightarrow  f_* \Theta_X   \rightarrow 
(\Theta_Y)  \rightarrow $$
$$  \rightarrow (\V^{\vee} \otimes W) \rightarrow  \mathcal{R}^1 f_*
\Theta_X  
\rightarrow  (\Theta_Y)\otimes \V^{\vee}  \rightarrow $$
$$ \rightarrow (\Lambda^2(\V^{\vee})\otimes W ) \rightarrow  \mathcal{R}^2 f_*
\Theta_X  
\rightarrow  (\Theta_Y)\otimes \Lambda^2(\V^{\vee}) \ra ,$$
where we used that  $\mathcal{R}^2 f_*(\hol_X) \cong \Lambda^2(\V^{\vee}).$

Notice that the coboundaries are given by cup product with the Kodaira-Spencer
class. Thus, in case where we have a holomorphic bundle, all rows
are exact.

\begin{oss}
In general, by a small variation of a theorem of E. Horikawa concerning the
deformations of holomorphic maps, namely Theorem 4.9 of
\cite{cat91}, we obtain in particular that, under the assumption $H^0
((\Theta_Y)\otimes V^{\vee} ) =0$,  we have a smooth morphism $Def (f)
\rightarrow Def(X)$ . 
\end{oss}

Let's now assume that we have a holomorphic principal torus bundle, thus
$W$ and $V$ are trivial holomorphic bundles on $Y$.

\begin{prop}
Let $f : X \ra Y$ be a holomorphic principal torus bundle and assume that
$H^0(\Theta_Y) =0$. Then every small deformation of $X$ is a holomorphic
principal torus bundle over some small deformation of $Y$.
\end{prop}

\Proof
Since $\V$ is trivial, by the previous remark, every  small deformation of $X$
is induced by a deformation of the map $ f : X \ra Y$.

That is, if we consider the respective Kuranishi families, we have a holomorphic
maps of $Def(f)$ to $Def(X)$, $Def(Y)$, such that the first is onto and smooth.

On the other hand, $H^1(f^* W)$ is the subspace $T^1(X|Y)$ in Flenner's notation
( cf. \cite{flen79}), kernel of the tangent map $ Def(f) \ra Def(Y)$, thus we
infer by the smoothness of $ Def(f) \ra Def(X)$, and since 
we have an exact sequence 
$$  0 \ra H^1(f^* W)  \ra H^1(\Theta_X ) \ra H^1(\Theta_Y) \ra  \dots ,$$ 
that we have indeed an isomorphism $ Def(f) \cong Def(X)$.

Then we have a morphism $Def(f) \cong Def(X) \ra Def(Y)$ and therefore, since any
deformation of a torus is a torus,
we get that any deformation of $X$ is a torus bundle over a deformation of $Y$.

However, the local system $\HH$ is trivial, and therefore it remains trivial also
after deformation. Whence, it follows that also the bundles $V$, $W$ remain
trivial, since they are both generated by $d$ global sections, by virtue
of (3) and (5). 

Therefore the Jacobian bundles remain trivial, so all the fibres are isomorphic
and every small deformation is a holomorphic bundle. 

Since  the classifying class in $H^1(Y,\M)$ is locally constant, we
finally infer that we get only principal holomorphic bundles.

\qed

\begin{oss}
 The cohomology group $ H^1(f^* W)$, by the Leray
spectral sequence, has a filtration whose graded quotient maps as a subspace of
$$ H^1(W) \oplus  H^0(Y, \RR ^1f_* \hol_X \otimes W) = H^1(\hol_Y^d) \oplus 
H^0(Y,\hol_Y^{d} \otimes \hol_Y^{d})= H^1(\hol_Y^d) \oplus \C^{d \times d}.$$

The meaning of the above splitting is explained by the exact sequence (9):
the left hand side stands for the deformation of the given bundle with fixed
fibre, the right hand side stands for the local deformation of the torus $T$.
\end{oss}

A typical case where the previous result applies is the one where $Y$ is a curve
of genus $g \geq 2$.

\begin{teo}
Let $f : X \ra Y$ be a holomorphic principal torus bundle over a curve $Y$ of
genus $g \geq 2$, and with fibre a $d$-dimensional complex torus $T$. Then the
Kuranishi family of deformations of
$X$  is smooth of dimension $3g - 3 + dg + d^2$, it has a smooth fibration onto
$ Def(Y) \times Def (T)$, its fibres over a point $(Y', T')$ are given by the
deformations parametrized by $ H^1 ( {\hol_{Y'}}^d)$ and inducing  all the
holomorphic principal torus bundles over $Y'$ with fibre $T'$ and given
topological class in
$H^2(Y', \Z^{2d})$.

\end{teo}

\Proof
We can use proposition 3.4 to construct, fixed a multidegree 
$m:= (m_1, .. m_d)$, a family $\FF _{g,m}$ with smooth base,
parametrizing quotients $X = L / N$, where  L is a $(\C^*)^d$-bundle with
multidegree $m$ over a smooth curve $Y$.

Here, the parameter $N$ varies in a smooth connected
$d^2$-dimensional parameter space for complex tori $(\C^*)^d / N$, while
the pair $(Y,L)$ varies in the fibre product of the Picard bundles $
Pic^{m_i}$ over a universal family parametrizing all smooth curves of
genus $g$.

By the proof of our previous proposition, our smooth family
$\FF _{g,m}$ has bijective Kodaira-Spencer map, whence it is locally
isomorphic to the Kuranishi family $Def(X)$.

The other assertions follow then easily.

\qed

\begin{oss}
Looking more carefully, we see that the choice of a  multidegree 
$m:= (m_1, .. m_d)$ is not unique. Indeed, one has a well defined element
in $H^2 (Y, \Lambda)  = H^2( Y, \Z) \otimes \Lambda \cong  \Lambda$,
where the isomorphism is unique since a complex structure on $Y$ fixes the orientation.
Thus we have as invariant a vector in a lattice $\Lambda \cong \Z^{2d}$,
whose only invariant is the divisibility index $\mu$ ( we define $\mu = 0$ if the
class is zero). 

Whence, we can reduce ourselves to consider only families $\FF _{g,\mu}$.
\end{oss}

\section{Deformations in the large of torus bundles over curves }

As was already pointed out, the 
examples  of Blanchard, Calabi and Sommese (cf. especially \cite{somm75})
show  that a manifold diffeomorphic to a
product $C \times T$,
where
$C$ is a curve of genus $g\geq 2$, and $T$ is a complex torus,
need not be a  holomorphic torus bundle over a curve. The deformation 
theory of Blanchard-Calabi varieties and, more generally, of differential
torus bundles offers interesting questions (cf. \cite {cat02}).

But, for the more narrow class of varieties which are holomorphic torus
bundles over a curve of genus $\geq 2$ (this includes as a special case
the products $C\times T$, where
$T$ is a complex torus, which were shown in \cite {cat02} to be closed 
for taking limits)  we have a  result  quite similar to the one we
have for tori.

\begin{teo}
A  deformation in the large of a holomorphic principal torus bundle over
a curve $C$ of genus $g\geq 2$ with fibre a complex torus $T$ of
dimension $d$ is again a holomorphic principal torus  bundle
over a curve $C'$ of genus $g$.

This clearly follows from the more precise statement.
Fix an integer 
$g\geq 2$, and a  a multidegree 
$m:= (m_1, .. m_d) \in \Z^d$: then every compact complex manifold $X_0$
such that its Kuranishi family of
deformations
$\pi : \Xi \rightarrow \mathcal{B}$ enjoys the property that the set
$\mathcal{B}" : = \{b|\ X_b$ is isomorphic to  a  manifold in the family
$\FF _{g,m}$  considered in Theorem 4.4 $\}$ has
$0$ as a limit point.
Then $X_0$ is also isomorphic to a manifold in the family
$\FF _{g,m}$ .
\end{teo}

PROOF.
Observe first of all that $\mathcal{B}"$ is open in the Kuranishi family
$\mathcal{B}$,
by the property that the Kuranishi family induces a versal family in each
neighbouring point, and by theorem 4.4.

Whence, by lemma  $2.8$ we can limit ourselves to consider the situation
where $\mathcal{B}$ is a 1-dimensional disk, and there exists a sequence 
$t_{\nu} \ra 0$
such that $X_{t_{\nu}}$,
$\forall \nu$, is a principal holomorphic torus bundle over a
curve $C_{t_{\nu}}$.

{\bf CLAIM I}
The first important property is that, for such a holomorphic principal
bundle $X = X_{t_{\nu}}$ as above, the subspace
$f_{t_{\nu}}^*(H^1(C_{t_{\nu}},\C))  $ is a
 subspace $V$ of $H^1( X,\C)  $,
 of dimension $2g$, such that each isotropic subspace $U$ of $H^1( X,\C) 
$ (i.e., such that the image 
 of $\Lambda^2(U) \rightarrow H^2( X,\C)  $ is zero) is contained in $V$.
 
 Moreover, if the bundle is not topologically trivial, $V$ itself is an
isotropic subspace, whereas, if the bundle is topologically trivial,
then any maximal isotropic subspace has dimension $g$.

In both cases, if  an isotropic subspace $U$ has dimension $r$ and $U +
\bar{U}$ has dimension $2r$, then $ r \leq g$.

{\bf Proof of claim I.}
Assume that  $ U $ is isotropic.
Notice that we have an injection $U \subset  V
\oplus W$, where $W$ injects into $ H^1(F, \Q)$ by virtue of the Leray
spectral sequence yielding 
$$ H^1(X, \Q) / H^1(Y, \Q) = ker( H^0(Y, \RR ^1
f_*(\Q))
\ra H^2(Y, \Q) ) .$$ Observe that $ H^0(Y, \RR ^1 f_*(\Q)) = H^1 (F,
\Q)$ since we have a principal bundle, and moreover that the above
homomorphism in the spectral sequence is determined by the cohomology class 
$\epsilon \in H^2 (Y, \Lambda)$.

In particular, the above map is  trivial if and only if the bundle is
topologically trivial.

In this case follows by the K\"unneth formula that the wedge product
yields an injection $ (W \otimes V ) \oplus (\Lambda^2 W) \ra H^2( X,\C)$ .

Whence, if  we have two non proportional cohomology
classes with trivial wedge product in cohomology, $(v_1 + w_1)
\wedge (v_2 + w_2) = 0$ and we assume w.l.o.g. $w_2 \neq 0$, then first of
all $ \exists \ c \in \C  \ s.t. \ w_1 = c w_2$.

Then also $c \neq 0$, $v_1 = c v_2$, a contradiction.

Assume now the bundle not to be topologically trivial: then 
$H^0(Y, \RR ^1 f_*(\C)) \ra H^2(Y, \C) ) $ is non trivial and $H^2(Y, \C)$
maps to zero in the cohomology of $X$, thus $V$ is an isotropic subspace.

However, then the Leray spectral sequence at least guarantees that
$ H^1 (Y,H^0(Y, \RR ^1 f_*(\C)) =  V \otimes H^1 (F, \C)$ is 
a direct summand in $H^2(X, \C)$.

Since moreover $W$ is killed by the linear form $\epsilon :
H^1 (F,\Q) \ra H^2(Y, \C) \cong \C$, $(\Lambda^2 W)$ embeds into the quotient 
$H^2(X, \C) / V \otimes H^1 (F, \C)$.

Thus again we have an injection $ (W \otimes V ) \oplus (\Lambda^2 W) \ra
H^2( X,\C)$, and the same argument as in the trivial case applies.

{\bf STEP II}. There is a morphism $F: \Xi \rightarrow \mathcal{C}$, where 
$\mathcal{C} \rightarrow \mathcal{B}$ is a smooth family of curves of genus
$g$.

{\bf Proof of step II.} 

We use for this purpose ideas related to the classical Castelnuovo - de
Franchis Theorem and to the isotropic subspace theorem (cf. \cite{cat91}).

From the differentiable triviality of our family $\Xi \rightarrow \Delta$
and by Claim I follows that we have a uniquely determined  subspace $V$ of
the cohomology of $X_0$, that we may freely identify to the previously
considered subspace $V$, for each
$X_{t_{\nu}}$.

Now, for each ${t_{\nu}}$, we have a decomposition $V = U_{t_{\nu}} \oplus
\overline {U_{t_{\nu}}}$, where $ U_{t_{\nu}} =
f^*(H^0(\Omega^1_{C_{t_{\nu}}}))  .$

By compactness of the Grassmann variety and by the  weak 1-Hodge
decomposition in the limit, the above decomposition also holds for $X_0$,
and $U_0$ is an isotropic subspace in  $ H^0( d \hol_{X_0} )$.
The Castelnuovo de Franchis theorem applies, and we get a holomorphic map
$f_0$ to a curve $Y_0$ of genus $\geq g$. But this genus must equal $g$,
since we noticed in the statement of Claim I that if  $U_0$ isotropic, and
$ U_0 + \bar{U_0} = U_0 \oplus \bar{U_0} $, then $dim (U_0) \leq g$.

Consider the vector bundle $f_* \Omega^1_{\Xi | \Delta}$.
As in the proof of Lemma 2.4 we infer the existence of a vector
(sub-)bundle $\HHH$ such that its stalk at $t$ yields a subspace of 
$ H^0( d \hol_{X_t} )$.

We have a map of the complex vector bundle $\HHH \bigoplus \bar{\HHH}$
to the trivial vector bundle $ H^1( X_t , \C ) = V \oplus W'$.

Since the construction of the Kuranishi family can be done by fixing the
underlying real analytic structure, we may assume w.l.o.g. that this map
of complex vector bundles has real analytic coefficients.

Now,$\forall t$, by Step I, $ Ker ( \HHH_t \ra W'_t)$ has dimension $
\leq g$, and dimension equal to $g$ for $ t = 0$.

Whence, this kernel provides a subspace $U_t$ $\forall t$ in a
neighbourhood of $0$. The corresponding map to the Grassmann manifold is
holomorphic in a non empty open set and real analytic, and is therefore
holomorphic $\forall t$ in a neighbourhood of $0$. 

We can now put the maps $f_t$ together by choosing a basis of 
$f_* (U_t)$, and integrating these holomorphic 1-forms: we get the
desired morphism $F: \Xi \rightarrow \mathcal{C}$ to the desired family of
curves of genus $g$.

 {\bf FINAL STEP}. 
 
We  have produced  a morphism $F: \Xi \rightarrow \mathcal{Y}$
where $\mathcal{Y} \rightarrow \mathcal{B}$ is a family of curves.

Moreover, the fibres of $f_0$ are deformations in the large of complex
tori, whence, by theorem 2.1, these fibres are complex tori of dimension
$d$. Consider, $\forall t_{\nu} $ the Kodaira -Spencer map
of
$Y_{t_{\nu}}$: it is identically zero in a neighbourhood of
$t_{\nu} $, thus it is identically zero in a neighbourhood of $0$.

Whence, $f_0$ is also a holomorphic torus bundle, and it is also
principal because the cohomology group $ H^1(Y_t, \M)$ is a subspace
of the cohomology group $ H^1(Y_t, Aut (\HH))$, and as such the
cohomology classes whose triviality ensure that a bundle is principal
are locally constant in $t$.

\hfill Q.E.D. \\ 

\begin{oss}
Consider the family $\FF _{g,m}$ of principal holomorphic torus bundles:
it is an interesting question to decide which manifolds in the given
family are K\"ahler manifolds, and which are not.
\end{oss}

\section{Holomorphic torus bundles over tori.}

In this section we shall consider the case where we have a holomorphic
principal bundle $ f:  X \ra Y$ with base a complex torus 
$ Y = V / \Gamma$ of dimension $m$, and fibre a complex torus
$T = U / \Lambda$ of dimension $d$.

In this case, our manifold $X$ is a $K (\pi, 1)$ and, as we already saw,
its fundamental group $ \pi_1 (X) := \Pi $ is a central extension
$$ (10') 1 \ra \Lambda \ra \Pi \ra \Gamma \ra 1$$
where the extension of $\Gamma$ by $\Lambda$ is classified
by  a cohomology class
$\epsilon \in H^2(Y, \Lambda )$.

Tensoring the above exact sequence with $\R$, we obtain an exact sequence
of Lie Groups
$$ (10') 1 \ra \Lambda \otimes \R \ra \Pi \otimes \R \ra \Gamma 
\otimes \R \ra 1$$ such that as a differentiable manifold our $X$
is the quotient $$  M:=  \Pi \otimes \R /  \Pi .$$

We want to give a holomorphic family of complex structures on $M$ such
that all holomorphic
principal bundles $ f:  X '\ra Y'$ with extension class isomorphic
to $\epsilon$ occur in this family.

This construction is completely analogous to the construction of the
standard family of complex tori, parametrized by an open set  $\mathcal
T_n$ in the complex
Grassmann Manifold $Gr(n,2n)$.

To this purpose, we must explain the analogue of the First Riemann
bilinear Relation in this context.

In abstract terms, it is just derived from the exact cohomology 
sequence
$$(9')   H^1 (Y, \hol_Y^d) \cong H^1 (Y, \HHH(U)_Y) \ra H^1(\HHH(T)_Y)
\ra ^c \ra H^2(Y, \Lambda ) \ra H^2 (Y, \HHH(U)_Y) $$
 telling that the class $\epsilon$ maps to zero in $H^2 (Y, \HHH(U)_Y) $.
This condition can be in more classical terms interpreted as follows:

\begin{oss}
{\bf( First Riemann Relation for  prin. hol. Torus Bundles)} 

There exists an
alternating bilinear map $A : \Gamma \times \Gamma \ra \Lambda$,
(representing the cohomology class $\epsilon$ uniquely)
 such that, viewing $A$ as a real element of 
 $$\Lambda^2(\Gamma 
\otimes \R )^{\vee} \otimes (\Lambda \otimes \R )=  \Lambda^2(V \oplus
\bar{V})^{\vee} \otimes (U \oplus \bar{U}),$$
its component in $  \Lambda^2(
\bar{V})^{\vee} \otimes (U )$ is zero.
\end{oss}

\Proof
By Dolbeault's theorem, the second cohomology group $H^2 (Y, \HHH(U)_Y) $
of the sheaf of holomorphic funtions with values in the complex vector
space $U$, since $Y$ is a torus, is isomorphic to the space of alternating
complex antilinear functions on $V \times V$ with values in $U$.
It is easy to verify that the coboundary operator corresponds to 
the projection of the cohomology group 
$$H^2(Y, \Lambda  \otimes \C) \ra
H^2 (Y,
\HHH(U)_Y) \cong   \Lambda^2 (\bar{V})^{\vee} \otimes (U ).$$
\qed

\begin{df}
Given $A$ as above, we define $\mathcal
T \mathcal B_A$ as the subset of the product of 
Grassmann Manifolds $Gr(m,2m)\times Gr(d,2d)$ defined by
$\mathcal
T \mathcal B_A : = \{ (V,U) | \ {\rm the \ component \ of\ A  \ in \ }
(\Lambda^2(
\bar{V})^{\vee} \otimes (U ))  \ is \ = 0 \}.$
\end{df}

\begin{oss}
Observe that  $\mathcal
T \mathcal B_A$ is a complex analytic variety, of codimension at most
$ d \ \frac {m (m-1)} {2} $.
\end{oss}
\Proof
We may in fact choose a basis for $\Gamma$,
resp.
for $\Lambda$  so that $A$ is then represented by a tensor $A^k_{i,j}$.
Let as usual $V$ correspond to a $2m \times m$ matrix $V$ ( thus the 
matrix $( V , \bar{V})$ yields the identity of  $\Gamma \otimes \C$ with
 arrival basis the chosen one and with initial basis $v_1,  \dots  v_m,
  \bar{v_1},  \dots \bar{v_m}.$
 
 Our desired condition is that, 
 $$ \forall h \neq \ell ,\ w_{h, \ell}:=  \ [ \Sigma_k (\Sigma_{i,j}  \
v_{i,h}  A^k_{i,j} v_{j, \ell} ) e_k]  \in U.$$

But this condition (observe moreover that $w_{h, \ell} = - w_{ \ell, h}$)
means that each vector $w_{h, \ell}$ is linearly dependent upon $u_1, 
\dots  u_d$. Locally in the Grassmannian $Gr(d,2d)$ the condition is
given by $d$ polynomial equations.
\qed

\begin{df}
The standard (Appell-Humbert) family of torus bundles parametrized
by $\mathcal T \mathcal B_A$ is the family of principal holomorphic torus
bundles $X_{V,U}$ on $Y : = V/ \Gamma$ and with fibre $T: = U /
\Lambda$ determined by the cocycle in
$H^1 (\Gamma ,\HHH(T)_Y) $ which is gotten by taking $ f_{\gamma} (z)$
which is the class 
${\rm mod} \ (\Lambda)$ of 
$ F_{\gamma} (z) : = - A (z, \gamma)  \ , \forall z \in V.$ 

\end{df}

\begin{oss}
Observe that we may write $A$ as $ B + \bar{B}$, where 
$$B \in 
 [ \Lambda^2 ({V})^{\vee} \otimes (U )]  \oplus [( V \otimes
\bar{V})^{\vee} 
\otimes (U )].$$
For later use, we write $B = B' + B"$, with $B' \in 
 [ \Lambda^2 ({V})^{\vee} \otimes (U )] ,\ B" \in \ [( V \otimes
\bar{V})^{\vee} 
\otimes (U )],$ and will say that $B"$ is the Hermitian component of $A$,
and $B'$ is the complex component of $A$.

Clearly,
\begin{itemize}
\item
$ A (z, \gamma) = B (z, \gamma)  \ , \forall z \in V,$
thus $ F_{\gamma} (z)$ is complex linear in $z$ with values in $U$.
\item
$ F_{\gamma} (z)$ is a cocycle with values in $T: = U /
\Lambda$ since 
$$ F_{\gamma_1 +\gamma_2 } (z) -  F_{\gamma_1} (z+ \gamma_2) - 
F_{\gamma_2} (z) = - F_{\gamma_1} (\gamma_2) = A (\gamma_1, \gamma_2) \in
\Lambda.$$
\end{itemize}
\end{oss}

\begin{df}
The complete Appell-Humbert family of torus bundles parametrized
by $\mathcal T' \mathcal B_A$ is the family of principal holomorphic torus
bundles $X_{V,U, \phi}$ on $Y : = V/ \Gamma$ and with fibre $T: = U /
\Lambda$ determined by the cocycle in
$H^1 (\Gamma ,\HHH(T)_Y) $ which is gotten by taking  the sum of 
$f_{\gamma} (z)$ with any cocycle $\phi \in H^1 (\Gamma ,\hol_Y^d) = H^1
(\Gamma ,\HHH(U)_Y).$
\end{df}

\begin{oss}
Observe that any $\phi \in H^1 (\Gamma ,\hol_Y^d) = H^1
(\Gamma ,\HHH(U)_Y)$ is represented by a cocycle with constant values in
$U$. Therefore the cocycle $f_{\gamma} (z) + \phi$ is a linear cocycle
(it is a polynomial of degree at most $1$ in $z$).
\end{oss}

\begin{teo}
Any  holomorphic principal torus bundle with extension class isomorphic
to $\epsilon  \in H^2 (\Gamma, \Lambda)$ occurs in the complete
Appell-Humbert family  $\mathcal T' \mathcal B_A$.
\end{teo}

\Proof
 Using sheaf cohomology, the proof is a standard consequence of the
exactness of sequence $(9)'$. 

\qed

\begin{oss}
The previous theorem is  an extension to the torus bundle case of
the classical theorem of linearization of the system of exponents for
holomorphic line bundles on complex tori ( cf. \cite{sieg}, pages 49-62
for an elementary proof, and \cite{c-c}, page 38, for a simpler elementary
proof).

The above elementary proofs do not work verbatim here, since \cite{sieg}
uses integrals which cannot be taken with values in $T$, and \cite{c-c}
uses the existence of a maximal  isotropic sublattice of $\Gamma$
( this no longer exists if $A$ has values in a lattice $\Lambda$ of
arbitrary rank).

Let us now indicate the modifications required in order to obtain an
elementary proof of theorem 6.8.

\end{oss}

\Proof n. II 

We have a cocycle $f_{\gamma}(z)$ with values in $T= U/\Lambda$, and
we can lift it to a holomorphic function $F_{\gamma}(z)$ defined
on the complex vector space $V$ and with values in $U$.

The cocycle condition tells us that 

$$ F_{\gamma_1 +\gamma_2 } (z) -  F_{\gamma_1} (z+ \gamma_2) - 
F_{\gamma_2} (z) = - F_{\gamma_1} (\gamma_2) = A (\gamma_1, \gamma_2) \in
\Lambda,$$ therefore $d F_{\gamma}(z)$ is a cocycle with values in 
$ H^0 (V, d (\hol_V^d) ) \subset H^0 (V, \hol_V^d)$.

By the classical linearization theorem there is a holomorphic $1$-form
$g$ such that 
$$ (II-1) \ d F_{\gamma}(z) = g (z + \gamma) - g (z) + L_{\gamma}
(z),$$ where  $L_{\gamma} (z)$ is linear (polynomial coefficients of degree
at most $1$).

Define $\psi (z)  := dg (z)$: by (II-1) follows that $\psi(z + \gamma) -
\psi (z)$ has constant coefficients, thus the coefficients of $\psi$
are linear (since their derivatives are $\Lambda$-periodic).

Since moreover $ d \psi = d^2 (g) = 0$, it follows that there is a
holomorphic 1- form $Q$ with quadratic coefficients such that 
$ d Q = \psi $.

Since $ d (g - Q) = 0$, there is a holomorphic function $\Phi (z)$ such
that $ d \Phi = g - Q$.

Now, $F_{\gamma}$ is cohomologous to $P_{\gamma }: = F_{\gamma }- \Phi ( z
+\gamma) +\Phi (z) $. Now, $ d P_{\gamma } = L_{\gamma}
(z) + Q(z) - Q (z + \gamma)$, thereby proving

{\bf Claim 1} The cocycle $f_{\gamma}$ is cohomologous to a cocycle
represented by  $P_{\gamma }$, whose coefficients are polynomials of
degree at most $3$.

We end the proof via the following

{\bf Claim 2}  The cocycle $f_{\gamma}$ is cohomologous to a cocycle
represented by  $G_{\gamma }$, whose coefficients are polynomials of
degree at most $1$.

Here, the proof runs as classically, since the relation (valid for every
cocycle)

$$P_{\gamma_1 } ( z + \gamma_2) - P_{\gamma_1 } ( z )
= P_{\gamma_2 } ( z + \gamma_1) -P_{\gamma_2 } ( z ) $$
shows that, if we choose a basis of $\C^m$ where $\Gamma$ has basis
$ e_1, \dots  e_m, \tau_1, \dots \tau_m$, then if $P_{e_i}$ is linear
$\forall i=1, .
\dots m$, it follows then that any $P_{\gamma }$ is linear.

To linearize $P_{e_i}, \forall i = 1, \dots m,$ one can use the methods
of solving linear difference equations as in \cite{sieg}, pages 54-58.

\qed

The explicit standard normal form for cocycles is very useful
to calculate several holomorphic cohomology groups of $X$:

\begin{teo}
The cokernel of $ 0 \ra H^0 (\Omega^1_Y )\ra H^0 (\Omega^1_X )$ is the
subspace of ${U}^{\vee}$ which annihilates the image of the Hermitian
part of $A$, i.e.,  of the component $B"$ in 
$[( V \otimes \bar{V})^{\vee}  \otimes (U )].$

It follows in particular that $X$ is parallelizable if and only if the
Hermitian part of $A$ is zero.
\end{teo}

\Proof
We have the exact sequence 
$$ 0 \ra H^0 (\Omega^1_Y )\ra H^0 (\Omega^1_X ) \ra H^0(f_*
\Omega^1_{X|Y}) \cong U^{\vee}$$
thus we immediately get that our cokernel is a subspace of the complex
vector space $U^{\vee}$.

To see which subspace do we get, let us take a holomorphic 1-form $\omega$ and lift
it back to the product $V \times T \cong \C^m \times T$, with coordinates
$(z,u)$.

Hence we write $\omega = \alpha + \beta = \Sigma_i \alpha_i (z,u) d
{z}_i + \Sigma_j \beta_j (z,u) d {u}_j $.  Since the coefficients are
holomorphic functions of $u$, and by the previous exact sequence we infer
 that $\beta$ has constant coefficients, and that the
coefficients of $\alpha$ depend only on $z$, we can write:
$\omega = 
\Sigma_i \alpha_i (z) d
{z}_i + \Sigma_j \beta_j  d {u}_j $. We also know that if $\beta = 0$ then
the coefficients $\alpha_i (z)$ are constant.

The condition that $\omega$ is a pull-back from $X$ is equivalent to
$\gamma^* (\omega) = \omega \ \forall \gamma \in \Gamma$.

Since $ \gamma (z,u) = ( z + \gamma, u + f_{\gamma} (z))$,  our
condition reads out as 
$$ (**)\ \alpha (z + \gamma) - \alpha (z) = - \beta \circ d f_{\gamma} (z)
: = - 
\Sigma _k \beta_k A^k (dz,\gamma).
$$ 
(**) implies that the coefficients of $\alpha$ may be assumed w.l.o.g. to
be linear homogenous functions of $z$, thus there is a complex bilinear
function
$B^*$ such that $ \alpha = B^* ( dz, z)$. Comparing again equation (**), we
infer $  B^* ( dz, \gamma)  = - \beta \circ A (dz,\gamma)$.

Since however we have $ \beta \circ A (dz,\gamma) =  \beta \circ B
(dz,\gamma)$ the above equation may hold iff $  B^* ( dz, \gamma)  = -
\beta \circ B' (dz,\gamma)$ and $\beta \circ B" (z,\gamma) = 0 \ \forall
\gamma, \ \forall z. $ Thus there exists an $\omega$ with vertical part $=
\beta$ if and only if $\beta$ annihilates the image of $B"$, as claimed.

\qed
\begin{cor}
The space $H^0 ( d\hol_X)$ of closed holomorphic 1-forms on $X$ contains 
the pull- back of $H^0 ( \Omega^1_Y)$ with cokernel the subspace $U^*$
of $U^{\vee}$ which annihilates the image of $A$, i.e.,  $U^* = \{ \beta | 
\beta
\circ A (z,\gamma)= 0 \
\forall
\gamma, \ \forall z\} .$
\end{cor}

\Proof
We know that there exists a holomorphic 1-form $\omega$ with vertical part
$\beta$ iff $\beta$ annihilates the Hermitian part $B"$ of $A$, and
then we may assume $\omega = 
\beta \circ B' (dz,z) + \Sigma_j \beta_j  d {u}_j .$ 

Let us calculate the differential of $\omega$: we get $ d \omega = 
\beta \circ B' (dz,d z) $. Since however $B'$ is antisymmetric, 
we get that $ d \omega = 0$ if and only if  $\beta \circ B' (z',z)=
0 \ \forall z, z' \in V$. Putting the two conditions together we obtain the
desired conclusion that $ \beta$ should annihilate the total image of $A$.

\qed

\begin{teo}
The cokernel of $ 0 \ra H^1 (\hol_Y )\ra H^1 (\hol_X )$ is the subspace
of $\bar{U}^{\vee}$ which annihilates the image of the anti-complex component  
 of $A$, i.e., of the conjugate of the component $B'$ in 
$[( \Lambda^2 V )^{\vee}  \otimes (U )].$
\end{teo}

\Proof
We calculate $H^1 (\hol_X )$ through the Leray spectral sequence
for the map $f$, yielding the exact sequence
$$ (*)\  0 \ra H^1 (\hol_Y )\ra H^1 (\hol_X ) \ra H^0( \RR^1 f_* \hol_X )
\ra H^2 (\hol_Y ).$$

We will interpret the above as Dolbeault- cohomology groups, keeping in
mind that $\RR^1 f_* \hol_X $ is a trivial holomorphic bundle of rank 
$d= dim (U)$.

Let $\eta$ be a $\bar{\partial}$-closed $(0,1)$-form, which we represent 
through its pull back to $ V \times T$, with coordinates $(z,u)$.

Hence we write $\eta = \alpha + \beta = \Sigma_i \alpha_i (z,u) d
\bar{z}_i + \Sigma_j \beta_j (z,u) d \bar{u}_j $.

The meaning of (*) is that is that $\eta$ is $\bar{\partial}$-cohomologous
to a form such that the functions $\beta_j$ are constant.
Let us then assume this to be the case, and let us impose the two
conditions $\bar{\partial} (\eta) = 0$, and that $\eta$ be a pull back
from $X$. 

The first implies $\bar{\partial}_u (\alpha) = 0$, thus the functions
$\alpha_i (z,u) $ are holomorphic in $u$, hence they only depend upon $z$.

Moreover then $\bar{\partial} (\eta) = 0$ is equivalent to
$\bar{\partial}_z (\alpha)(z) = 0$, and therefore is equivalent to the
existence of a function
$\phi(z)$ on $\C^n$ such that  $\bar{\partial}_z  (\phi(z)) = (\alpha)(z).$

Then, the condition that $\eta$ be a pull back means that $\eta$
is invariant under the action of $\Gamma$ such that
$$  z \ra z + \gamma , \  u \ra  u + f_{\gamma}(z) ,$$
whence it writes out as 
$ \alpha (z + \gamma) - \alpha (z) = -  \Sigma_k
\beta_k  \overline { d f_{\gamma}(z)_k }$, or, equivalently,
$$ (***) \ \bar{\partial}  \phi(z + \gamma) - \bar{\partial}  \phi(z) =
-  \Sigma_k
\beta_k  A^k  (\gamma , d\bar{z}) .$$

As before, we can write $A$ as $A' + A"$ where, $B"$ being the Hermitian part
of the tensor, $A" = B" + \bar{B"}$.

Assume first that $ \beta \circ B' (\gamma ,
\bar{z}):=  \Sigma_k \beta_k  {B'}^k  (\gamma ,
\bar{z}) $ is identically zero: then $ \beta \circ B" (\gamma ,
\bar{z}) = \beta \circ A (\gamma ,
\bar{z})$ and we simply choose $ \phi: = \beta \circ B" (z ,
\bar{z})$ which clearly satisfies equation (***).

Then $ \eta : = \bar{\partial} ( \phi (z)) + \Sigma_j \beta_j  d \bar{u}_j $
is a $\bar{\partial} $-closed 1-form on $X$ with vertical part $\beta =
\Sigma_j \beta_j  d \bar{u}_j$.

Conversely, the meaning of the exact sequence (*) is as follows:

given any  vertical part $\beta = \Sigma_j \beta_j  d \bar{u}_j$,
we first find a differential $(0,1)$-form on $X$

$\eta = \alpha + \beta = \Sigma_i \alpha_i (z) d
\bar{z}_i + \Sigma_j \beta_j  d \bar{u}_j $ with the given vertical part.

Then, the vertical part comes from $H^1 (\hol_X)$ only if its image
in $H^2 (\hol_Y)$, provided by $\bar{\partial} (\eta)$, is zero.

For this purpose we have to solve the equation

$$    \alpha_i (z + \gamma) - \alpha_i (z) =
-  \Sigma_k
\beta_k  A^k  (\gamma , \bar{e_i})= - \beta \circ A (\gamma , \bar{e_i}) .$$

But an easy solution is provided by setting

$$    \alpha_i (z) =
-  \Sigma_k
\beta_k  {B'}^k  (\bar{z} , \bar{e_i}) -  \Sigma_k
\beta_k  {B"}^k  (z , \bar{e_i}) = - \beta \circ B' (\bar{z},\bar{e_i}) 
- \beta \circ B" (z,\bar{e_i}).$$

But now, $\bar{\partial}  (\eta) =\bar{\partial} \alpha =  - \beta \circ B'
(d \bar{z},d \bar{z})$ which is a (0-2)-form with constant coefficients, and 
is therefore equal to zero if and only if all of its coefficients are zero,
thus iff $ - \beta \circ B'
( \bar{v},d \bar{zw}) = 0  \ \forall v,w \in V$, which proves our assertion. 

\qed

\begin{oss}
The previous theorem allows us to write explicitly the basis of the complete
Appel-Humbert family.  Locally on $Gr(d, 2d)$ we  may assume after a permutation
of coordinates  that our  subspace $U$ is the space of vectors $\{(u', u") | u" =
U^* (u')\}$, for some $(d \times d)$-matrix $U^*$.

Then our equations are: we consider the matrices $(V,U^*)$ and the
vectors $v
\in \C^m$, $\beta
\in \C^d$ such that, in the notation of Remark 6.3,
 $ \forall h \neq \ell ,\ w_{h, \ell}:=  \ [ (\Sigma_{i,j}  \
v_{i,h}  A_{i,j} v_{j, \ell} ) ]    $ satisfies
$$ w_{h, \ell}" = U^* (w_{h, \ell}'), \ \ \beta (w_{h, \ell}') = 0.   $$ 

Let us consider the easiest possible case, where $m=2, d=1$: then our parameter
variety is a pull back of the variety in $\C^4$ with equations
$$ w" = u^* w', \ \ \beta w' = 0,   $$
a product of $\C$ with $\{ (\beta, w') \in \C^2 | \ \beta w' = 0 \}   $,
 a reducible variety.
\end{oss}

In order to analyse the problem of describing the small deformations of
principal holomorphic torus bundles over tori we need to calculate the
cohomology groups of the tangent sheaf $\Theta_X$. This can be accomplished by the
above multilinear algebra methods.

For ease of calculations recall that the canonical bundle $\Omega^n_X $ is trivial,
therefore, by Serre duality, it is sufficient to determine the cohomology groups  
$$ H^{n-i} (\Omega^1_X) \cong  H^i (\Theta_X)^{\vee} .$$ 

\begin{teo}
$ H^{n-i} (\Omega^1_{X|Y})$ fits into a short exact sequence 
 $$ 0 \ra coker \beta_{n-i-1} \ra H^{n-i} (\Omega^1_{X|Y}) \ra ker \beta_{n-i} \ra
0,$$  where $\beta_i : U^{\vee} \otimes H^i (\hol_X) \ra V^{\vee} \otimes H^{i+1}
(\hol_X)$ is given by cup product and contraction with $ B" \in  [(  \bar{V})
\otimes V )^{\vee} \otimes (U )]$.
\end{teo}
\Proof
Recall the exact sequence (6)
$$(6) \  0 \ra f^* \Omega^1_Y \ra  \Omega^1_X \ra f^* W^{\vee} \ra 0$$
where $W$ is a trivial bundle, and the extension class is a pull back of the
extension class of the sequence
$$(8) \  0 \ra  \Omega^1_Y \ra f_* \Omega^1_X \ra W^{\vee} = f_* \Omega^1_{X|Y}
\ra 0 ,$$
which coincides with the Hermitian component $B"$ of $A$, as it is easy to verify
( observe that $ B" \in  [(  \bar{V}) \otimes V )^{\vee} \otimes (U )] \cong
H^1_{\bar{\partial}} ( V^{\vee}  \otimes U  \otimes \hol_Y) \cong H^1( W
\otimes\Omega^1_Y)  \subset H^1( f^* W
\otimes f^*\Omega^1_Y) $).

Therefore, in the exact cohomology sequence of (6), the coboundary operator is given
by cup product with $B"$, whence we can write this  exact sequence as 
$$ \dots H^{n-i} (\Omega^1_X) \ra H^{n-i} (\Omega^1_{X|Y}) = H^{n-i} ( U^{\vee}
\otimes
\hol_X) 
\ra H^{n-i+1} (f^*\Omega^1_Y) =  H^{n-i+1} ( V^{\vee} \otimes
\hol_X) 
$$ 
whence the desired assertion follows.

\qed

\begin{cor}
$H^i (\Theta_X)$ fits into a short exact sequence
$$ 0 \ra coker b_{i-1} \ra H^i (\Theta_X) \ra ker b_{i} \ra
0,$$ where  $b_i : V \otimes H^{i} (\hol_X) \ra U \otimes H^{i+1}
(\hol_X)$ is given by cup product and contraction with $ B" \in  [(  \bar{V})
\otimes V )^{\vee} \otimes (U )]$.
\end{cor}
\Proof
Serre duality.
\qed

For completeness we state without proof a more general result than what we need

\begin{prop}
The Leray spectral sequence for the sheaf $\hol_X$ and for the map $f$ yields a
spectral sequence  which degenerates at the $E_3$ level and with $E_2$ term  $ = (H^i
(\RR^j f_*
\hol_X ) , d_2 )$, where $d_2 : (H^i (\RR^j f_* \hol_X )= H^i ( \Lambda^j (V)) = 
\Lambda^i (\bar{V}^{\vee}) \otimes \Lambda^j (\bar{U}^{\vee}) \ra \Lambda^{i+2}
(\bar{V}^{\vee}) \otimes \Lambda^{j-1} (\bar{U}^{\vee})$ is provided by cup product
 and contraction with $\bar{B'} \in \Lambda^2 (\bar{V}^{\vee}) \otimes (\bar{U})$. 
\end{prop}

We observe now some deformation theoretic consequences:

\begin{teo}
Let $ f:  X \ra Y$ be a holomorphic principal bundle with base a complex torus 
$ Y = V / \Gamma$ of dimension $m$, and fibre an elliptic curve $T$
($T = U / \Lambda$ has dimension $1$).

Assume moreover $ \pi_1 (X) := \Pi $ to be a nontrivial  central extension
$$ (10') 1 \ra \Lambda \ra \Pi \ra \Gamma \ra 1$$
classified
by  a cohomology class
$\epsilon  \neq 0 \ \in H^2(Y, \Lambda )$
whose associated bilinear form $A$ has an image of dimension $=2$.

1) Then every limit of  manifolds in the complete Appell-Humbert family is again a
holomorphic principal bundle $ f' :  X' \ra Y'$ with fibre an elliptic curve $T'$,
and thus occurs in the complete Appell-Humbert family (6.6 and 6.8).

2) Every small deformation of such a manifold $X$ is induced by the complete
Appell-Humbert in the case where $h^1(\hol_X) = m$.
\end{teo}

\Proof
By corollary 6.11 we get that $H^0 ( d\hol_X)$, the space  of closed holomorphic
1-forms on
$X$ equals  the pull- back of $H^0 ( \Omega^1_Y) = H^0 ( d\hol_Y)$ since the
antisymmetric form $A$ representing $\epsilon  $ is non zero.

Likewise, again the Leray spectral sequence for cohomology shows that 
$ H^1 (X, \C)$ contains the pull back of $ H^1 (Y, \C)$ with cokernel 
equal to the subspace of $\Lambda^{\vee}$ annihilated by the image of
$A$. This image is the whole $\Lambda$ , whence $b_1(X) = b_1(Y)$ and
it follows immediately that
 $ f:  X \ra Y$ is the Albanese map of $X$. 
 
 {\bf  Proof of assertion 1}: 
 By lemma 2.4 and cor. 2.5 any limit $X'$  
in  a 1-parameter family $X_t$  of such holomorphic bundles has a
surjective Albanese map 
$ f' :  X' \ra Y'$ onto a complex torus 
$ Y' = V' / \Gamma$ of dimension $m$.

Since any deformation in the large of an elliptic curve is an elliptic curve,
the general fibre of $f'$ is an elliptic curve.

We are done if we show that $f'$ is a submersion, since an elliptic fibration
without singular fibres is a holomorphic bundle. Moreover, since we have that the
small deformations of  $ f' :  X' \ra Y'$ contain principal holomorphic
bundles, then the same thing holds for $ f' :  X' \ra Y'$.

That $f'$ is a submersion follows by purity of branch locus (cf. \cite {mum78}): we
would have then a ramification divisor $R_0$ for the central fibre, but
$R_{t_{\nu}} =
\emptyset$  for the special fibres, yielding the same contradiction as in the
proof of theorem 2.7.
 
{\bf  Proof of assertion 2}: 

It suffices to show that every small deformation $X'$ of $X$ has an Albanese
map  $ f' :  X' \ra Y'$ onto a complex torus 
$ Y' = V' / \Gamma$ of dimension $m$.

However, by semicontinuity, $h^1(\hol_{X'}) \leq m$, on the other hand we may
apply lemma 2.3 to conclude that $X'$ has a very good Albanese map, and again
by purity of branch locus we conclude that we have a submersion.

\qed

\begin{oss}
We are interested more generally to see whether the complete Appell Humbert family,
which does not have a smooth base, but is pure of dimension $m^2 + m$, has
surjective Kodaira-Spencer map (this is a first step towards the proof of
its versality). It is necessary for this purpose to calculate the cohomology
group
$H^1 (\Theta_X)$ in an explicit way.

We write, according to theorem 6.12,  $H^1 (\hol_X ) = (\bar{V}^{\vee}) \oplus ker
(B')$, where we see
$B'$ as a linear map $ (B'): (\bar{U}^{\vee})  \ra \Lambda^2 (\bar{V}^{\vee})$.

We have then , according to corollary 6.14, an exact sequence 
$$ V \ra U \otimes (\bar{V}^{\vee} \oplus ker
B') \ra H^1 (\Theta_X) \ra  V \otimes (\bar{V}^{\vee} \oplus ker
B')   \ra  U \otimes ( \Lambda^2 (\bar{V}^{\vee}) \oplus 
(\bar{V}^{\vee} \otimes (ker B'))$$
where the first and last map are given by contraction with $B"$.

We first  analyse when holds the
desired dimensional estimate $ dim H^1 (\Theta_X) = m^2 + m $.

Case 1 : $B"$ is zero, whence $X$ is parallelizable, thus $ dim H^1 (\Theta_X) =
(m+1) h^1 (\hol_X) $. But in this case $B'$ is non zero, whence $h^1 (\hol_X)= m  $
and we are done. 

Case 2: $B" \neq 0$  and $B' \neq 0$, then   $h^1 (\hol_X)= m  $ and thus 
follows $ dim H^1 (\Theta_X) \leq  (m+1) h^1 (\hol_X) = m (m+1)$.

We omit the not difficult verification that in both cases the 
Kodaira Spencer map is  surjective. 

Case 3. $B' = 0$ , and $B" \neq 0$.
In this case $A$ = $B" + \bar{B"}$, whence  we observe  that the image of $B"$
has dimension $=1$. 

We leave aside this case, where plays a role also the rank of the antisymmetric
matrix $B"$. 
\end{oss}

We finally describe the most well known  example of such elliptic bundles
over 2-dimensional tori, namely, the so called {\bf Iwasawa manifold} (cf.
\cite{k-m71}.

One lets  $N $ (biholomorphic to $  \C^3$) be the unipotent group of $3
\times 3$ upper trangular matrices with all eigenvalues equal to $1$.

$N$ contains the discrete cocompact subgroup $\Pi$ which is the subgroup of the
matrices with entries in the subring $\Lambda \subset \C$,  $\Lambda : =
\Z [i]$.

Obviously, $X:= N / \Pi$ is a fibre bundle $ f : X \ra \C^2/ \Lambda^2$, where
$f$ is induced by the homomorphism $F : N \ra \C^2$ given by the coordinates
$z_{12}, z_{23}$. 

$X$ is parallelizable being a quotient of a complex Lie group by a discrete
subgroup.

Moreover, in our case, the alternating form $A$ is gotten as the
antisymmetrization of the product map $\C \times \C \ra \C$ , which
induces  $\Lambda \times \Lambda \ra \Lambda$.

Thus in this case there is no Hermitian part, confirming the parallelizability, and
then $h^1 (\hol_X) = 2$ , and $dim H^1 (\Theta_X) = 6$.

That's why the Kuranishi family of $X$ is smooth of dimension $6$, and  gets small
deformations which are not parallelizable (cf.\cite{nak75}).

\bigskip

\section{ Blanchard -Calabi Manifolds}

The  Sommese-Blanchard examples (\cite{bla53}, \cite{bla56},
\cite{somm75}) provide non K\"ahler complex structures
$X$ on manifolds diffeomorphic to a product $C \times T$, where $C$ is a
compact complex curve and $T$ is a 2-dimensional complex torus.

In fact, in these examples, the projection $X \rightarrow C$ is holomorphic 
and all the fibres are 2-dimensional complex tori.

\begin{oss}

 Start with a Sommese-Blanchard 3-fold with trivial canonical
bundle and deform the curve $C$ and the line bundle $L$ in such a way that
the canonical bundle becomes the pull back of a non torsion element of
$Pic^0(C)$: then this is the famous example showing that the Kodaira
dimension is not deformation invariant for non K\"ahler manifolds ( cf.
\cite{ue80}).
\end{oss}

Also Calabi (\cite{cal58})   showed that there are  non K\"ahlerian 
complex structures  on a product $C\times T$ .

We briefly sketch  (cf.\cite{cat02} for precise results in the case
where the base manifold $Y$ is a curve  and X has dimension 3) how to
vastly generalize these constructions.

\begin{df} (BLANCHARD-CALABI Jacobian MANIFOLDS)

Let $Y$ be a compact complex manifold, and let $W$ be a rank $d$
holomorphic vector bundle admitting $2d$ holomorphic sections
$\sigma_1,\sigma_2, ..  ,\sigma_{2d}$ which are everywhere
$\R$-linearly independent.

Then the quotient $X$ of $W$ by the
$\Z^{2d}$-action acting fibrewise by translations : $ w \rightarrow  w +
\Sigma_{i=1,..2d}\ n_i \sigma_i$ is a complex manifold diffeomorphic to 
the differentiable manifold $Y \times T$,  where $T$ is a $d$- dimensional
complex torus. $X$  will be called a {\bf Jacobian Blanchard-Calabi
manifold}.
\end{df}

\begin{oss}
The canonical divisor of $X$ equals $K_X =
\pi^* (K_Y - det W) $, where $\pi : X \rightarrow Y$ is the canonical
projection. Moreover, if $Y$ is K\"ahler, $X$ is K\"ahler  only if the
bundle $W$ is trivial (i.e., iff $X$ is  a holomorphic product $Y \times
T$).

Indeed, if $d=2$, one
has $h^0 (\Omega^1_X) = h^0 (\Omega^1_Y) $  unless $W$ is trivial,
as follows from the dual of the sequence (5)
$$ 0 \rightarrow V \rightarrow (\hol_Y)^{2d}
\rightarrow  W \rightarrow 0.$$
\end{oss}

\Proof

 If $X$ is K\"ahler, then, since the first Betti number of $X$ equals
$2d+ 2 h^0 (\Omega^1_Y)$, then it must hold that $h^0 (\Omega^1_X) = h^0
(\Omega^1_Y) + d$, in particular
$h^0 (W^{\vee}) \geq d$ .

But $(W^{\vee})$ is a subbundle of a trivial bundle of rank $2d$,
whence $d$ linearly independent sections of $(W^{\vee})$ yield 
a composition $(\hol_Y)^d \rightarrow (W^{\vee}) \rightarrow
(\hol_Y)^{2d}   $ whose image is a trivial bundle of rank $\leq d$.

Indeed the rank must be $d$, else the $d$ sections would not be
$\C$-linearly independent.

Assume now that $h^0 (W^{\vee}) =r$: then the corresponding image of the
composition $(\hol_Y)^r \rightarrow (W^{\vee}) \rightarrow
(\hol_Y)^{2d}   $ would be a trivial bundle of rank 
$r $, thus  $(W^{\vee})$ has a
trivial summand $I_r$ of rank $r$ . 
Then  we have a direct sum $(W^{\vee}) = I_r \oplus Q$  and dually $W$ is
a direct sum $I_r \oplus  Q^{\vee}$.  We use now the hypothesis that
there are $2d-2r$ holomorphic sections  of
$W$ which are $\R$-independent at any point: it follows then that
if $d=2$ and $0 < r < d=2$, then $Q^{\vee}$ has rank 1 and admits 
holomorphic sections  which are everywhere
$\R$-independent, thus it is trivial.

\qed \\

\begin{df}
Given a Jacobian Blanchard-Calabi manifold  $ \pi : X \rightarrow Y$,
any $X$-principal homogeneous space $\pi' :Z \rightarrow Y$ will be called
a Blanchard-Calabi manifold.
\end{df}

\begin{oss}
Given a Jacobian Blanchard-Calabi manifold $X \ra Y$, the associated $X$
- principal homogeneous spaces are classified by $H^1(Y, W)$.

Thus, the main existence problem is the one concerning the Jacobian
Blanchard-Calabi manifolds.
\end{oss}

\begin{prop}
Let $G$ be the  $d^2$- dimensional projective Grassmann variety
$G \PP(d-1,2d-1)$ of $d-1$-dimensional projective subspaces in
$\PP^{2d-1}$. The datum of Jacobian Blanchard-Calabi manifold $ f : X \ra
Y$ is equivalent to the datum of a "totally non real" holomorphic  map
$ h : Y \ra G$, i.e., such that the corresponding ruled variety
$R_h$ given by the union of the $d-1$-dimensional projective subspaces
$h(y)$ has no real points.
\end{prop}

\Proof

Let $G$ be the Grassmann variety $G \PP(d-1,2d-1) = Gr (d, 2d)$, so $G$
parametrizes the
 $d$-dimensional vector  subspaces in $\C^{2d}$:
observe that the datum of an exact sequence 
$$ 0 \rightarrow V \rightarrow (\hol_C)^{2d} \rightarrow  W \rightarrow
0$$ is equivalent to the datum of a holomorphic map $h: Y\rightarrow G$,
since for any such $f$ we let $V,W$ be the respective pull backs of the
universal subbundle $U$ and  of the quotient bundle $Q$.

The condition that the standard $2d$ sections are $\R$-linearly
independent means that there is no $y \in Y$ and no real vector 
$v \in \R^{2d}$ such that $v$ belongs to the subspace corresponding
to $h(y)$. 
\qed

If  $h$ is  constant,   $W$ is  trivial and we have a product. 

Else, we  assume for simplicity that $h$ is generically finite
onto its image, and we consider the deformation of the holomorphic map
$h$.

We have the following exact sequence:
$$ 0 \rightarrow \Theta_Y \rightarrow (f)^* \Theta_G  \rightarrow  N_h
\rightarrow 0 ,$$
where $N_h$, the normal sheaf  of the morphism $h$, governs the
deformation theory of the morphism $h$, in the sense that the tangent
space to Def(h) is the space $H^0(N_h)$, while the obstructions lie in
$H^1(N_h)$.

By virtue of the fact that $\Theta_G  = Hom(U, Q)$, and of the
cohomology sequence associated to the aboce exact sequence, we get 
$ 0 \rightarrow H^0(\Theta_Y ) \rightarrow H^0(V^{\vee} \otimes W) 
\rightarrow  H^0(N_h) \rightarrow H^1(\Theta_Y ) \rightarrow H^1(V^{\vee}
\otimes W) \rightarrow  H^1(N_f) 
\rightarrow 0 $, and we conclude that the deformations of the map are
unobstructed provided $H^1(V^{\vee} \otimes W) = 0$.

In \cite{cat02} we saw that this holds, in the special case where $W = L
\oplus L$, if the degree
$d$ of $L$  satisfies $d \geq g$. 
We showed in this way that if $d \geq g$ the dimension of the space of
deformations of the map
$f$ is given by $ 3g - 3 + 4 h^0(2L) = 4 d + 1-g $, and this number
clearly tends to infinity together with $ d = deg(L)$.

Whence we got (loc. cit.) the following
\begin{cor}
The space of complex structures on the product of a curve $C$
with a four
dimensional real torus has unbounded dimension. 
\end{cor}

As a corollary of our theorem 4.1 we  also obtain

\begin{cor}
The space of complex structures on the product of a curve $C$
of genus $g \geq 2$ with a four
dimensional real torus contains manifolds which are not deformation
equivalent to each other, namely, Blanchard-Calabi 3-folds which are not
K\"ahler and holomorphic principal bundles in the family $\FF _{g,0}$.
\end{cor}

\bigskip

{\bf Acknowledgements.}

I wish to thank  Fedya Bogomolov and Hubert Flenner for some useful
remarks,   Paola Frediani and Alessandra Sarti for  some
stimulating discussions during the Wintersemester 1997-98 in
G\"ottingen.

\vfill

\noindent
{\bf Author's address:}

\bigskip

\noindent 
Prof. Fabrizio Catanese\\
Lehrstuhl Mathematik VIII\\
Universit\"at Bayreuth\\
 D-95440, BAYREUTH, Germany

e-mail: Fabrizio.Catanese@uni-bayreuth.de


\bigskip


  
\begin{thebibliography}{99}



\bibitem[BPV84]{bpv84} 

W. Barth, C. Peters and A. Van de Ven,
``{\em Compact Complex Surfaces}'', {\bf Ergebnisse der Mathematik und
ihrer Grenzgebiete, 3.F, B. 4, Springer-Verlag}, 1984. 


\bibitem[Bea78]{bea78}
A. Beauville, 
``{\em Surfaces algebriques complexes}'',
{\bf Asterisque  54} Soc. Math. Fr. (1978)

\bibitem[Bla53]{bla53}
A. Blanchard, ``{\em Recherche de structures analytiques complexes 
sur certaines vari\'et\'es}'', {\bf C.R. Acad.Sci., Paris 238} (1953), 
657-659.

\bibitem[Bla56]{bla56}
A. Blanchard, ``{\em Sur les vari\'et\'es
analytiques complexes}'', {\bf Ann.Sci. Ec. Norm. Super., III Ser., 73}
(1956), 157-202.


\bibitem[Cal58]{cal58}
E. Calabi, 
``{\em Construction and properties of some 6-dimensional almost complex
manifolds}, {\bf Trans. Amer. Math.Soc. 87  } (1958), 407-438.

\bibitem[Cat84]{cat84}
F. Catanese, 
``{\em On the Moduli Spaces of Surfaces of General Type}'',
{\bf J. Diff. Geom 19} (1984) 483--515.


\bibitem[C-C91]{c-c}
F. Capocasa - F. Catanese, 
``{\em Periodic meromorphic functions.}''
{\bf Acta Math. 166} (1991) 27-68.

\bibitem[Cat91]{cat91}
F. Catanese, 
``{\em Moduli and classification of irregular K\"ahler manifolds
 (and algebraic varieties) with Albanese general type fibrations.
 Appendix by Arnaud Beauville.}''
{\bf Inv. Math. 104} (1991) 263-289; Appendix 289 .


\bibitem[Cat95]{cat95}
F. Catanese, 
``{\em Compact
complex manifolds bimeromorphic to tori.}'' Proc. of the Conf. "Abelian
Varieties", Egloffstein
1993, {\bf  De Gruyter }(1995), 55-62.

\bibitem[Cat00]{cat00}
F. Catanese, 
``{\em Fibred surfaces, varieties isogenous to a product and 
related moduli spaces}''
{\bf Amer. Jour. Math. 122}  (2000)  1-44.

\bibitem[Cat96]{cat96}
F. Catanese,
{\em Generalized Kummer
surfaces and differentiable structures on Noether -Horikawa surfaces, 
I.} in 'Manifolds and
Geometry, Pisa 1993', {\bf Symposia Mathematica XXXVI, Cambridge University
Press} (1996), 132- 177.

\bibitem[Cat02]{cat02}
F. Catanese,
{\em Deformation types of  real and complex manifolds},
 in
'Contemporary Trends in Algebraic Geometry and Algebraic
Topology', Proc. of the Chen-Chow Memorial Conference, Nankai Tracts in Math.
5, Chern, Fu, Hain editors, World Scientific (2002), 193-236.



\bibitem[C-F02]{c-f02}
F. Catanese, P. Frediani,
``{\em Real hyperelliptic surfaces and the orbifold fundamental group}''
  math.AG/0012003,  
to appear in Journal of the Inst. Math. Jussieu {\bf 2} (1) (2003),
1-65.


\bibitem[Che58]{che58} 
S.S. Chern, 
``{\em Complex manifolds}'',
{\bf Publ. Mat. Univ.  Recife} (1958).


\bibitem[FG65]{fg65}
W. Fischer, H. Grauert, 
``{\em Lokal-triviale Familien kompakter komplexer Mannigfaltigkeiten }''
{ \bf Nachr. Akad. Wiss. G\"ottingen, II. Math.-Phys. Kl. 1965} (1965), 89-94 . 

\bibitem[Flen79]{flen79}
H. Flenner,
``{\em  \"Uber Deformationen holomorpher Abbildungen},
{\bf Osnabr\"ucker Schriften zur Mathematik }(1979), I- VII+ 1-142.


\bibitem[GR84]{gr84}
 H. Grauert, R. Remmert, 
``{\em Coherent analytic sheaves}''
{ \bf Grundlehren der Mathematischen Wissenschaften, 265,
 Springer-Verlag, XVIII}  (1984).

\bibitem [Hir62]{hir62}
H. Hironaka,
``{\em An example of a non-K\"ahlerian complex-analytic deformation of
K\"ahlerian complex structures}'', {\bf Ann. Math. (2) 75,}
(1962), 190-208.   
 



\bibitem [Ko63]{ko63} 
K. Kodaira ,''{\em On compact analytic
surfaces II-III}'', {\bf Ann.  Math.  77} (1963), 563-626, 
{\bf Ann.  Math.  78} (1963), 1-40.


\bibitem [Ko64]{ko64} 
K. Kodaira ,''{\em On the structure of complex analytic
surfaces I}'', {\bf Amer. J. Math.  86} (1964), 751-798 .

\bibitem [Ko68]{ko68} 
K. Kodaira ,''{\em On the structure of complex analytic
surfaces IV}'', {\bf Amer. J. Math.  90} (1968), 1048-1066 .

\bibitem [K-M71]{k-m71}
K. Kodaira, J. Morrow ,``{\em Complex manifolds}'' {\bf Holt,
Rinehart and Winston},  New York-Montreal, Que.-London (1971).

\bibitem[K-S58]{k-s58}
K. Kodaira , D. Spencer ''{\em On deformations of complex analytic
structures I-II}'', {\bf Ann. of Math.  67} (1958), 328-466 .





\bibitem[Maz63]{maz63} 
B. Mazur, ``{\em Differential topology from
the point of view of simple homotopy theory}'', {\bf Publ. Math. I.H.E.S.
21} (1963). 1-93, and correction ibidem 22 ( 1964) , 81-92 .  


\bibitem[Mos78]{mos78}
G. Mostow,
''{\em Strong rigidity of locally symmetric
spaces}'',  {\bf Annals of Math.Stud. 78, Princeton Univ. Press}(1978).

\bibitem[Mum70]{mum70}
D. Mumford, 
''{\em Abelian varieties}''
{\bf Tata Institute of Fundamental Research Studies in Mathematics, 
 Oxford University Press. VIII} (1970).

\bibitem[Mum78]{mum78}
D. Mumford, 
''{\em Some footnotes to the work of C. P. Ramanujam}''
{\bf C.P. Ramanujam. - A tribute. Collect. Publ. of C.P. Ramanujam 
and Pap. in his Mem., Tata Inst. fundam. Res., Stud. Math. 8}  (1978), 247-262.

\bibitem[Nak75]{nak75}
I. Nakamura, 
''{\em Complex parallelisable manifolds and their small deformations}''
{ \bf J. Differ. Geom. 10}, (1975), 85-112 .

\bibitem[Nak98]{nak98}
I. Nakamura, 
''{\em Global deformations of $\PP^2$-bundles over $\PP^1$}''
{ \bf J. Math. Kyoto Univ. 38, No.1} (1998), 29-54 .

\bibitem[Shaf61]{shaf61}
I.R. Shafarevich,''{\em Principal homogeneous spaces defined over a
function field}'' {\bf Steklov Math.Inst. 64 ( Transl. AMS vol. 37)}
(1961), 316-346 (85-115).

\bibitem[Sieg73]{sieg}
C.L. Siegel, ''{\em Topics in complex function theory, Vol. III}'' {\bf
Interscience Tracts in pure and appl. Math., n. 25, Wiley, N.Y.} (1973).

\bibitem[Somm75]{somm75}
A. J. Sommese ,``{\em Quaternionic manifolds}'', 
{\bf 
Math. Ann. 212}  (1975), 191-214.


  \bibitem [Su70]{su70} 
  T. Suwa, ``{\em On hyperelliptic surfaces}'', {\bf J. Fac. Sci.
Univ. Tokyo, 16} (1970), 469-476 .


 \bibitem [Su75]{su75} 
  T. Suwa, ``{\em Compact quotient spaces of $\C^2$ by affine transformation
groups}'',  {\bf J. Differ. Geom. 10}  (1975), 239-252.




\bibitem [Su77-I]{su77-I} 
  T. Suwa, ``{\em Compact quotients of $\C^3$ by affine transformation
groups, I.}'', {\bf Several complex Variables, Proc. Symp. Pure Math. 30, Part 1,
Williamstown 1975,} (1977), 293-295.




\bibitem [Su77-II]{su77-II} 
  T. Suwa, ``{\em Compact quotients of $\C^3$ by affine transformation
groups, II.}'', {\bf Complex Analysis and Algebraic Geometry, Iwanami
Shoten} (1977), 259-278 .

\bibitem [Ue75]{ue75} 
  K. Ueno, ``{\em
Classification theory of algebraic varieties and compact complex spaces.} 

{\bf Lecture Notes in Mathematics 439,} 
Springer-Verlag, XIX, 278 p.(1975).

\bibitem [Ue80]{ue80} 
  K. Ueno, ``{\em On three-dimensional compact complex manifolds with non
positive Kodaira dimension}'',  {\bf Proc. Japan Acad., vol.56, S.A.n.10}
(1980), 479-483.

 \bibitem [Ue82]{ue82} 
  K. Ueno, ``{\em Bimeromorphic Geometry of algebraic and analytic
threefolds}'', in 'C.I.M.E. Algebraic Threefolds, Varenna 1981'
{\bf Lecture
Notes in math. 947} (1982), 1-34 . 

\bibitem [Ue87]{ue87} 
  K. Ueno, ``{\em On compact analytic Threefolds with non trivial Albanese
Tori}'',  {\bf Math. Ann. 278} (1987), 41-70.




\end{thebibliography}
\end{document}